\documentclass{jimo}
\usepackage{amsmath}
\usepackage{subfig}
\usepackage{booktabs}
\usepackage{algorithmicx}
\usepackage{algpseudocode}
  \usepackage{paralist}
  \usepackage{graphics} 
  \usepackage{epsfig} 
 \usepackage[colorlinks=true]{hyperref}
\hypersetup{urlcolor=blue, citecolor=red}

     \def\DOI{}

\theoremstyle{definition}

\title[State Transition Algorithm]
      {State transition algorithm}

\author[Xiaojun Zhou, Chunhua Yang and Weihua Gui]{}


 \email{tiezhongyu2010@gmail.com}
 \email{ychh@csu.edu.cn}
 \email{gwh@mail.csu.edu.cn}

\thanks{Xiaojun Zhou is supported by the China Scholarship Council,
Chunhua Yang is supported by the National Science Found for Distinguished Young Scholars of China (Grant No. 61025015) and PCSIRT (IRT1044), who is the corresponding author of this paper, and Weihua Gui is supported by
the National Nature Science Foundation of China (Grant No. 61134006).
}

\begin{document}
\maketitle

\centerline{\scshape Xiaojun Zhou, Chunhua Yang and Weihua Gui}
\medskip
{\footnotesize
 \centerline{School of Information Science and Engineering}
   \centerline{Central South University, Changsha, 410083, China}
} 
\medskip
\bigskip
\begin{abstract}
In terms of the concepts of state and state transition, a new heuristic random search algorithm named state transition algorithm is proposed. For continuous function optimization problems, four special transformation operators called rotation, translation, expansion and axesion are designed. Adjusting measures of the transformations are mainly studied to keep the balance of exploration and exploitation. Convergence analysis is also discussed about the algorithm based on random search theory. In the meanwhile, to strengthen the search ability in high dimensional space, communication strategy is introduced into the basic algorithm and intermittent exchange is presented to prevent premature convergence. Finally, experiments are carried out for the algorithms. With 10 common benchmark unconstrained continuous functions used to test the performance, the results show that state transition algorithms are promising algorithms due to their good global search capability and convergence property when compared with some popular algorithms.
\end{abstract}

\section{Introduction}
The concept of state means to a situation which a material system maintains, and it is characterized by a group of physical qualities. The process of a system turning from a state to another is called state transition, which can be described by a state transition matrix. The idea of state transition was created by a Russian mathematician named Markov when he expected to represent a specific stochastic process (known as Markov process)\cite{bi1}. Not only in communication theory but also in modern control theory, state transition matrix is of great importance. For instance, in modern control theory, it can determine the stability of a system.\\
\indent In almost all branches of engineering, including system design, tactical planning, system analysis, process management and control, and model parameter adjustment, optimization techniques have found wide applications\cite{bi2}. Generally speaking, the methods used for solving such optimization problems can be classified into two categories: deterministic and stochastic, in which, stochastic methods are subdivided into evolutionary algorithms and metaheuristic algorithms. The traditional deterministic algorithms include Hooke-Jeeves pattern search\cite{bi3} and hill-climbing, evolutionary algorithms contain genetic algorithm (GA)\cite{bi4,bi5}, evolutionary programming, evolution strategies, and genetic programming, while metaheuristic algorithms consist of simulated annealing, particle swarm optimization (PSO)\cite{bi6,bi7}, differential evolution(DE) algorithm\cite{bi8}, etc. On the other hand, the commonly used numerical algorithms for engineering optimization problems can also be categorized into two classes: direct search methods and gradient-based methods. The direct search methods comprise the simplex search, Powell's conjugate direction method and random search, while the gradient-based methods include the Newton's (basic-, modified-, quasi-) and conjugate gradient methods\cite{bi9}. In the same time, hybrid methods, combining of deterministic and stochastic or direct search and gradient-based, are also proposed to  draw on each other's strengths\cite{bi10,bi11,bi12,bi13}.\\
\indent According to the No Free Lunch Theorem\cite{bi14}, no search algorithm is better than other algorithms on the space of all possible problems.
This paper introduces a new method for optimization of continuous nonlinear functions, which belongs to metaheuristic random search. Because of its foundation on state and state transition, the method is called state transition algorithm (STA)\cite{bi15,bi16}. The algorithm has roots in three main component methodologies. One is the random optimization theory, the others are population-based approach and then space transformation method. In this paper, it focuses on four operators named rotation, translation, expansion and axesion transformation as well as the communication strategy in state transition algorithm. Compared with some state-of-the-art optimization algorithms,  RCGA\cite{bi17}, CLPSO\cite{bi18}, and SaDE\cite{bi19}, which are
improved versions of GA, PSO and DE, the experimental results show that STAs are comparable and promising algorithms.
\section{The basic state transition algorithm}
Considering the following unconstrained optimization problem
\begin{equation}
\min_{x \in \Re ^{n}} f(x).
\end{equation}
\indent In a deterministic view, it usually adopts iterative method to solve the problem
\begin{equation}
x_{k+1}=x_{k} + a_{k}d_{k},
\end{equation}
where, $x_{k}$ is the \textit{k}th iteration point, $a_{k}$ is the \textit{k}th step size and $d_{k}$ is the \textit{k}th search direction.\\
\indent The common selection of a step is by exact line search or inexact line search. While the techniques of search direction include steepest descent method, conjugate gradient methods, Newton's methods, alternating directions, and conjugate direction methods\cite{bi20}.\\
\indent In a way, the iterative methods aim to search for a direction and a step in an iteration. Though these methods utilize the gradient information explicitly or implicitly, they have their inherent defects. For one thing, it is computationally difficult. For another, it only indicates the local information. In the view of global optimization, the direction of gradient is just a way standing for direction, and it has no substantial effect on searching for a global optimum. If the iterative method is concerned in a state and state transition way, then an iterative point can be regarded as a state, the process of searching for a direction and a step will equate to a state transition process, and through a state transition, a new state will be created.\\
\indent In the point of stochastic, it can also understand the evolutionary algorithms and metaheuristic algorithms in a state and state transition way. For example, genetic algorithm, its each individual of a generation can be considered as a state, and the updating process of using genetic operators such as selection, crossover and mutation can equate to state transition processes. In the same way, particle swarm optimization, the flock updating its velocity and position, and differential evolution, adding the difference vector of two randomly chosen vectors to a target vector, can also be regarded as state transition processes.\\
\indent In terms of the concept of state and state transition, a solution to a specific optimization problem can be described as a state, the operators of optimization algorithms can be considered as state transition, and the process to update current solution will become a state transition process.\\
\indent Through the above analysis and discussion, it defines the following form of state transition
\begin{equation}
\left \{ \begin{array}{ll}
x_{k+1}= A_{k}x_{k} + B_{k}u_{k}\\
y_{k+1}= f(x_{k+1})
\end{array} ,\right.
\end{equation}
where, $x_{k}$ stands for a state, corresponding to a solution to the optimization problem; then, $A_{k}$ and
$B_{k}$ are state transition matrixes, which can be regarded as operators of optimization algorithm;
$u_{k}$ is the function of state $x_{k}$ and historical states; while $f$ is the cost function or evaluation function.
\subsection{State transformation operators}
As a matter of fact, operators such as reflection, contraction, expansion and rotation are widely used in simplex optimization method\cite{bi21,bi22,bi23}, which is especially popular in the fields of chemistry, chemical engineering, and medicine. However, they always fail to lead to continued progress and are not applicable to a wide range of functions.\\
\indent In the theory of space and transformation, rotation matrices are only defined for two and three dimensional transformation. For example, the two dimensional rotation matrix is
\scriptsize
$\left[
\begin{array}{ll}
cos\theta & -sin\theta \\
sin\theta & cos \theta
\end{array}
\right]$. \\
\normalsize
\indent Using various types of space transformation for reference, in this paper, it defines the following four special state transformation operators to solve continuous function optimization problems.\\
(1) Rotation transformation
\begin{equation}
x_{k+1}=x_{k}+\alpha \frac{1}{n \|x_{k}\|_{2}} R_{r} x_{k},
\end{equation}
where, $x_{k}$ $\in$ $\Re^{n}$, $\alpha$ is a positive constant, called rotation factor; $R_{r}$ $\in$ $\Re^{n\times n}$, is random matrix with its entries obeying the uniform distribution in the range of [-1, 1] and $\|\cdot\|_{2}$ is 2-norm of vector or Euclidean norm. Then, it will prove that the rotation transformation has the function of searching in a hypersphere.
\begin{proof}
\begin{equation}
\begin{split}
\|x_{k+1}-x_{k}\|_{2} & = \|\alpha \frac{1}{n \|x_{k}\|_{2}} R_{r} x_{k} \|_{2}\\
                      & = \frac{\alpha}{n \|x_{k}\|_{2}} \|R_{r} x_{k}\|_{2}\\
                      & \leq \frac{\alpha}{n \|x_{k}\|_{2}} \|R_{r}\|_{m_\infty}\|x_{k}\|_{2} \leq \alpha
\end{split}
\end{equation}
\end{proof}
\noindent (2) Translation transformation
\begin{equation}
x_{k+1} = x_{k}+  \beta  R_{t}  \frac{x_{k}-x_{k-1}}{\|x_{k}-x_{k-1}\|_{2}},
\end{equation}
where, $\beta$ is a positive constant, called translation factor; $R_{t}$ $\in \Re$ is a random variable with its components obeying the uniform distribution in the range of [0,1]. It is obvious to find  the translation transformation has the function of searching along a line from $x_{k-1}$ to $x_{k}$ at the starting point $x_{k}$, with the maximum length of $\beta$.\\
(3) Expansion transformation
\begin{equation}
x_{k+1} = x_{k}+  \gamma  R_{e}x_{k},
\end{equation}
where, $\gamma$ is a positive constant, called expansion factor; $R_{e} \in \Re^{n \times n}$ is a random diagonal matrix with its elements obeying the Gaussian distribution (in this study, standard normal distribution). It is also obvious to find  the expansion transformation has the function of expanding the components in $x_{k}$ to the range of [-$\infty$, +$\infty$], searching in the whole space.\\
(4) Axesion transformation
\begin{equation}
x_{k+1} = x_{k}+  \delta  R_{a}x_{k},
\end{equation}
where, $\delta$ is a positive constant, called axesion factor; $R_{a}$ $\in \Re^{n \times n}$ is a random diagonal matrix with its entries obeying the Gaussian distribution and only one random position having nonzero value. The axesion transformation aims to search along the axes and strengthens single dimensional search.
\subsection{State transformation algorithm}
Before the state transition algorithm, it is necessary to introduce the basic random optimization\cite{bi24,bi25,bi26,bi27}. Considering the above unconstrained optimization problem, the procedure of the basic random optimization can be outlined in the following pseudocode.\\
\begin{algorithmic}[1]
\State Initialize feasible solution $x_{0}$, and set $k \gets 0$
\Repeat
   \State $k \gets k +1 $
   \State Generate a Gaussian random number vector $r$
   \State $x_{trail} \gets x_{k-1} + r$
   \If{$f(x_{trail})< f(x_{k-1})$}
   \State $x_{k} \gets x_{trail}$
   \Else
   \State {$x_{k} \gets x_{k-1}$}
   \EndIf
\Until{the specified termination criterion is met}
\end{algorithmic}
\quad \\
\indent For one thing, as a metaheuristic random method, the state transition algorithm is similar to the basic random optimization\cite{bi24}. The only difference is that a candidate solution set is generated by the four special operators, while a new trail is selected following the same way as that of the basic random optimization, which means that the ``greedy criterion'' is used in selecting the new state. By the way, a candidate solution set is created by some times of transformation. The times of the transformation or the size of the set is called search enforcement (\textit{SE}), and the translation operator is only performed when a better new trail is found.\\
\indent For another, as a stochastic algorithm\cite{bi28}, the dealing with dynamic balance between diversification (exploration of the solution space) and intensification (exploitation of the accumulated knowledge) is also significant in state transition algorithm. Due to their intrinsic properties, the rotation is chosen for exploitation, the expansion is for exploration, the translation is selected as to maintain equilibrium between them, and axesion is proposed to strength the single dimensional search.\\
\indent The main process of STA is shown in the pseudocode as follows
\begin{algorithmic}[1]
\Repeat
    \If{$\alpha < \alpha_{\min}$}
    \State {$\alpha \gets \alpha_{\max}$}
    \EndIf
    \State {Best $\gets$ expansion(funfcn,Best,SE,$\beta$,$\gamma$)}\Comment{expansion transformation}
    \State {Best $\gets$ rotation(funfcn,Best,SE,$\alpha$,$\beta$)}\Comment{rotation transformation}
    \State {Best $\gets$ axesion(funfcn,Best,SE,$\beta$,$\delta$)}\Comment{axesion transformation}
    \State {$\alpha \gets \frac{\alpha}{\textit{fc}}$}
\Until{the specified termination criterion is met}
\end{algorithmic}
\quad \\
\indent As for detailed explanations, expansion function in above pseudocode is given as follows for example
\begin{algorithmic}[1]
\State{oldBest $\gets$ Best}
\State{fBest $\gets$ feval(funfcn,oldBest)}
\State{State $\gets$ op\_expand(Best,SE,$\gamma$)}
\State{[newBest,fGBest] $\gets$ fitness(funfcn,State)}
\If{fGBest $<$ fBest}\Comment{greedy criterion}
    \State{fBest $\gets$ fGBest}
    \State{Best $\gets$ newBest}
    \State{State $\gets$ op\_translate(oldBest,Best,SE,$\beta$)}
    \State{[newBest,fGBest] $\gets$ fitness(funfcn,State)}
    \If{fGBest $<$ fBest}\Comment{greedy criterion}
        \State{fBest $\gets$ fGBest}
        \State{Best $\gets$ newBest}
    \EndIf
\EndIf
\end{algorithmic}
\subsection{Parameters analysis in STA}
In state transition algorithm, there are five important parameters, namely search enforcement(\textit{SE}), rotation factor $\alpha$, translation factor $\beta$, expansion factor $\gamma$ and axesion factor $\delta$. It is easy to understand that the larger the search enforcement, the higher the intensity of search, and vise versa. However, the larger search enforcement will cause larger computational complexity. In this paper, the search enforcement is recommended to use the same size as the dimension of the optimization problem.\\
\indent When \textit{SE} is constant, taking the exploration and exploitation into consideration, the strategy of adjusting parameters of the four operators is significant. To make the deeper exploitation, the smaller rotation factor is needed. Especially, the rotation factor will vary in a declining way from a positive constant till zero to gain a high precision solution. In the meanwhile, there are two schemes to regulate the $\alpha$. One is to adjust the parameter in an inner loop, namely, decreasing the rotation factor from a start constant to the end in the operation of rotation transformation\cite{bi15}. The other is to adjust the parameters in an outside loop, that is to say, decreasing the rotation factor according to the iterations. To balance the global search and local search timely, the latter scheme is adopted in the paper; however, the rotation factor is decreasing itself from a maximum value to a minimum value in an exponential way with base \textit{fc}, which is called lessening coefficient\cite{bi16}, as described in the pseudocode of STA. By the way, extra tests have testified the effectiveness of the scheme.\\
\indent As for the remained control parameters, for example, the larger the translation factor, the longer STA searches along a straight line. However, the magnitude of translation factor has great influence on exploitation and exploration. The large translation factor will facilitate the exploration, while the small translation factor benefits the exploitation. Similarly, the same phenomenon exists in the selection of expansion and axesion factors. Taking the complexity of adjusting strategies for these control parameters into consideration, we keep them fixed in current version of STA for simplicity.
\subsection{The convergence analysis of STA}
The convergence of stochastic optimization algorithms has been heatedly discussed. For instance, genetic algorithm was analyzed by means of homogeneous finite Markov chain\cite{bi29}, particle swarm optimization was studied to investigate particle trajectories in a discrete system view\cite{bi30}, and convergence property of differential evolution was also discussed in\cite{bi8}.\\
\indent As a metaheuristic random optimization algorithm, the convergence analysis of STA will follow the same way as random search methods. In fact, the probability of random search algorithm for finding global minimum being equal to 1 was stated by Solis and Wets\cite{bi31}. That is to say, the STA will satisfy the similar convergence performance of random optimization algorithm, and readers who are interested in convergence analysis are referred to their work for details.
\section{Communication strategy into state transition algorithm}
In a way, the basic state transition algorithm is individual-based, and an individual searches in its neighborhood. The difference between the basic state transition algorithm and other random optimizations is that the search space is normalized or specialized.\\
\indent The population-based approach is prevalent in metaheuristic algorithms, such as genetic algorithm, particle swarm optimization and differential evolution. Let's name the basic state transition algorithm STAI, the improved state transition algorithm based on population is called STAII with the number of states denoted as \textit{SN}. In the meanwhile, some communication strategies are necessary to manage the individuals for sharing information, which is important in population-based methods.
\subsection{Crossover operator}
Individual communication can be implemented in various ways, of which crossover operation is quite common, especially in genetic algorithms.\\
\indent Let $X_1$ and $X_2$ be individual components of current generation, $Y_1$ and $Y_2$ are the offspring components, some canonical crossover operators are displayed in the following.\\
(1)  Michalewicz's arithmetical crossover\cite{bi32}
\begin{equation}
\left \{ \begin{array}{ll}
Y_1= \alpha X_1 + (1-\alpha)X_2\\
Y_2= \alpha X_2 + (1-\alpha)X_1
\end{array} ,\right.
\end{equation}
where, $\alpha$ is either a constant or a variable whose value depends on the age of population.\\
(2) Wright's linear crossover\cite{bi33}
\begin{equation}
Y_1= 1.5 X_1 - 0.5 X_2, Y_2= -0.5 X_1 + 1.5 X_2, Y_3= (X_1+X_2)/2
\end{equation}
(3) Kalyanmoy Deb's simulated binary crossover\cite{bi34}
\begin{equation}
\left \{ \begin{array}{ll}
Y_1= 0.5[ (1-\beta)X_1 + (1+\beta)X_2]\\
Y_2= 0.5[ (1+\beta)X_1 + (1-\beta)X_2]
\end{array} ,\right.
\end{equation}
where, $\beta$ is a random variable, obeying the following probability distribution
\begin{equation}
\begin{cases}
p(\beta) = 0.5(\eta_c+1)\beta^{\eta_c}& 0\leq \beta \leq 1,\\
p(\beta) = 0.5(\eta_c+1)\frac{1}{\beta^{\eta_c + 2}}& \beta > 1,
\end{cases}
\end{equation}
here, $p(\cdot)$ is the probability density function, $\eta_c$ is the distribution index, which determine how well spread the children will be inherited from their parents.\\
(4) The proposed crossover
\begin{equation}
\left \{ \begin{array}{ll}
Y_1= \delta X_1 +(1-\delta) X_2,\\
Y_2= \eta X_1 + (1-\eta)X_2,
\end{array} \right.
\end{equation}
where, $\delta$ and $\eta$ are independent variables, which obey the 0-1 distribution.\\
\indent In this proposed crossover, crossover operation means for each component of a pair of individuals, components exchange or maintain their information completely.
\subsection{Intermittent exchange}
Different from other population based algorithms, all of the individuals in STAII are elites, and they develop themselves trough state transformation, which is referred as self learning. When the communication strategy is introduced, the individual can contact with each other to better develop themselves. However, it may bring about some disadvantageous effects. If the frequency of individual communication is too high, individuals are apt to imitate each other utterly, which will cause premature convergence. In this paper, intermittent exchange is proposed to solve the issue, that is to say, individual communication occurs at a certain frequency, where the frequency is named communication frequency (\textit{CF}).\\
\indent The communication strategy is adopted to share information among individuals, and it is regulated by communication frequency. If \textit{CF} is small enough, it will equate the situation without the communication strategy. When \textit{CF} is too large, individuals are easily trapped into imitating each other, causing the premature convergence, that is to say, a moderate \textit{CF} is appropriate. In this paper, we recommend to use the same magnitude as square root of the maximum iterations.\\
\indent When the exchange condition is satisfied, the proposed crossover operator will be performed. Each state will communicate with all of the other states,
to make sure that useful information is completely shared.
\subsection{The framework of STA with communication strategy}
Through the above discussion and analysis, by introducing in communication strategy, the pseudocode of the kernel of state transition algorithm
can be described as shown in the following
\begin{algorithmic}[1]
\Repeat
    \If{$\alpha < \alpha_{\min}$}
        \State {$\alpha \gets \alpha_{\max}$}
    \EndIf
    \State {State $\gets$ self\_learning(funfcn,State,SE,$\alpha$,$\beta$,$\gamma$,$\delta$)}\Comment{self learning}
    \State {$\alpha \gets \frac{\alpha}{\textit{fc}}$}
    \If{mod(iter,CF)==0}\Comment{intermittent exchange}
        \State {State $\gets$ communication(funfcn,State)}
    \EndIf
    \State {[Best,fBest] $\gets$ fitness(funfcn,State)}
\Until{the specified termination criterion is met}
\end{algorithmic}
\indent In the process of the algorithm, the $self\_learning$ means each state in the state set will be performed on four state transformation operators, while in $communication$ function, the intermittent exchange is adopted at intervals. The flowchart of the algorithm is outlined in Figure 1.
\begin{figure}[h!]
\centering
\includegraphics[width=8cm,height=10cm]{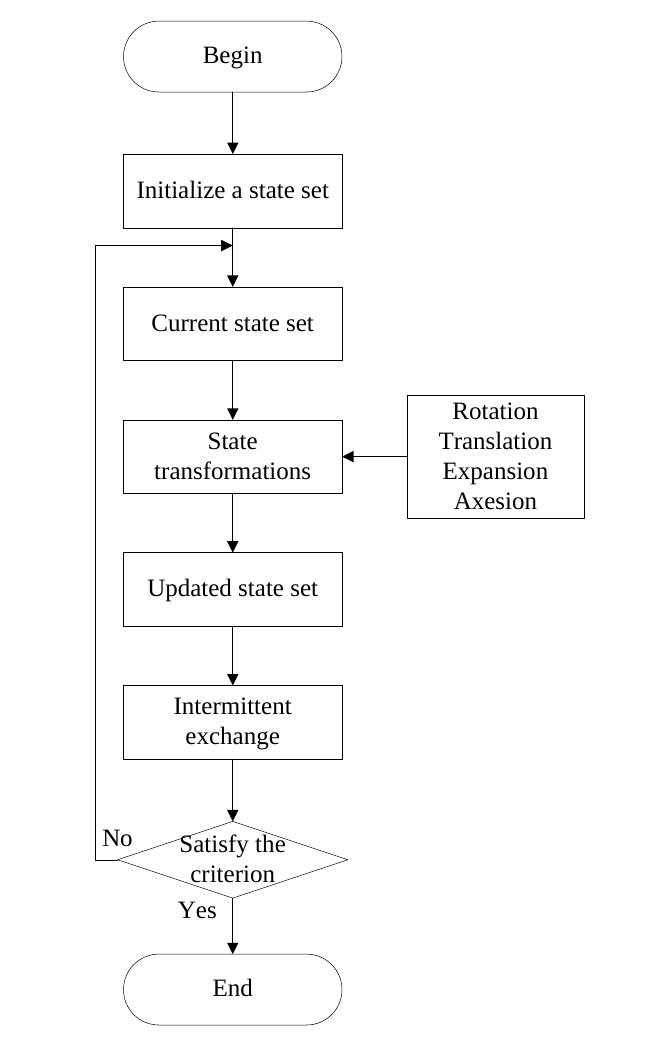}
\caption{the flowchart of STAII}
\end{figure}
\section{Experiments and results}
To compare the proposed state transition algorithm with previously mentioned RCGA, CLPSO and SaDE, two experiments are arranged. The first experiment is mainly for two dimensional functions, and the other focuses on ten dimensional functions test. In the same time, both STAI and STAII are carried out, for comparison with other algorithms as well as themselves.
\subsection{Test functions}
In order to test the performance of STA, ten common benchmark functions are selected for the experiment. Seven functions are multidimensional functions of various modals and the other three are two dimensional functions, which are listed in Table 1, while the landscapes of two dimensional functions are plotted in Figure 2.
\begin{table}[h!]
\caption{Benchmark functions for test in this paper}
\footnotesize 
\begin{tabular}{{cccc}}
\hline
\toprule[1pt]
Name of function & Function definition  & Range & $f_{min}$ \\
\hline
Spherical        & $~f_1= \sum_{i=1}^n x_{i}^{2} $  & [-100,100]    &  0\\
Rastrigin        & $~f_2=  \sum_{i=1}^n(x_{i}^{2}-10cos(2 \pi x_{i})+10)$ & [-5.12,5.12]  & 0\\
Griewank         & $~f_3=\frac{1}{4000} \sum_{i=1}^n x_{i}^{2}-  \prod_{i}^{n} cos|\frac{x_{i}}{\sqrt{i}}| + 1$  & [-600,600] & 0 \\
Rosenbrock       & $~f_4= \sum_{i=1}^n (100(x_{i+1}-x_{i}^2)^2 + (x_{i}-1)^2)$ & [-30,30] & 0 \\
Schewefel        & $~f_5= \sum_{i=1}^n [ -x_{i} sin(\sqrt{|x_{i}|})]$ & [-500,500] & -418.9829n \\
Ackley           &
$
\begin{array}{l}
f_{6}=20+e-20exp(-0.2\sqrt{\frac{1}{n} \sum_{i=1}^n x_{i}^2})\\
-exp(\frac{1}{n} \sum_{i=1}^n cos(2\pi x_{i}))
\end{array}
$
& [-32,32] & 0 \\
Michalewicz      & $~f_7= \sum_{i=1}^n sin(x_{i}) sin( \frac{i x_{i}^2}{\pi})^{20} $ & $[0,\pi]$ & - \\
Schaffer         & $~f_8= 0.5 + \frac{sin(\sqrt{x_{1}^2+x_{2}^2})^2-0.5}{(1+0.001(x_{1}^2+x_{2}^2))^2}$ & [-100,100] & 0 \\
Easom            &
$
\begin{array}{l}
f_{9}= -cos(x_{1})cos(x_{2})\times\\
exp(-(x_{1}-\pi)-(x_{2}-\pi))
\end{array}
$
& [-100,100] & -1 \\
Goldstein-Price  &
$
\begin{array}{l}
f_{10}= [1+(x_{1}+x_{2}+1)^2(19-14x_{1}+3x_{1}^2\\
-14x_{2}+6x_{1}x_{2}+3x_{2}^2)]\times [30+(2x_{1}-3x_{2})^2\\
(18-32x_{1}+12x_{1}^2+48x_{2}-36x_{1}x_{2}+27x_{2}^2)]
\end{array}
$
& [-2,2] & 3 \\
\bottomrule[1pt]
\hline
\end{tabular}
\end{table}
\begin{figure}[!htbp]
\centering
\subfloat[$f_1$]{\includegraphics[width=4cm,height=4cm]{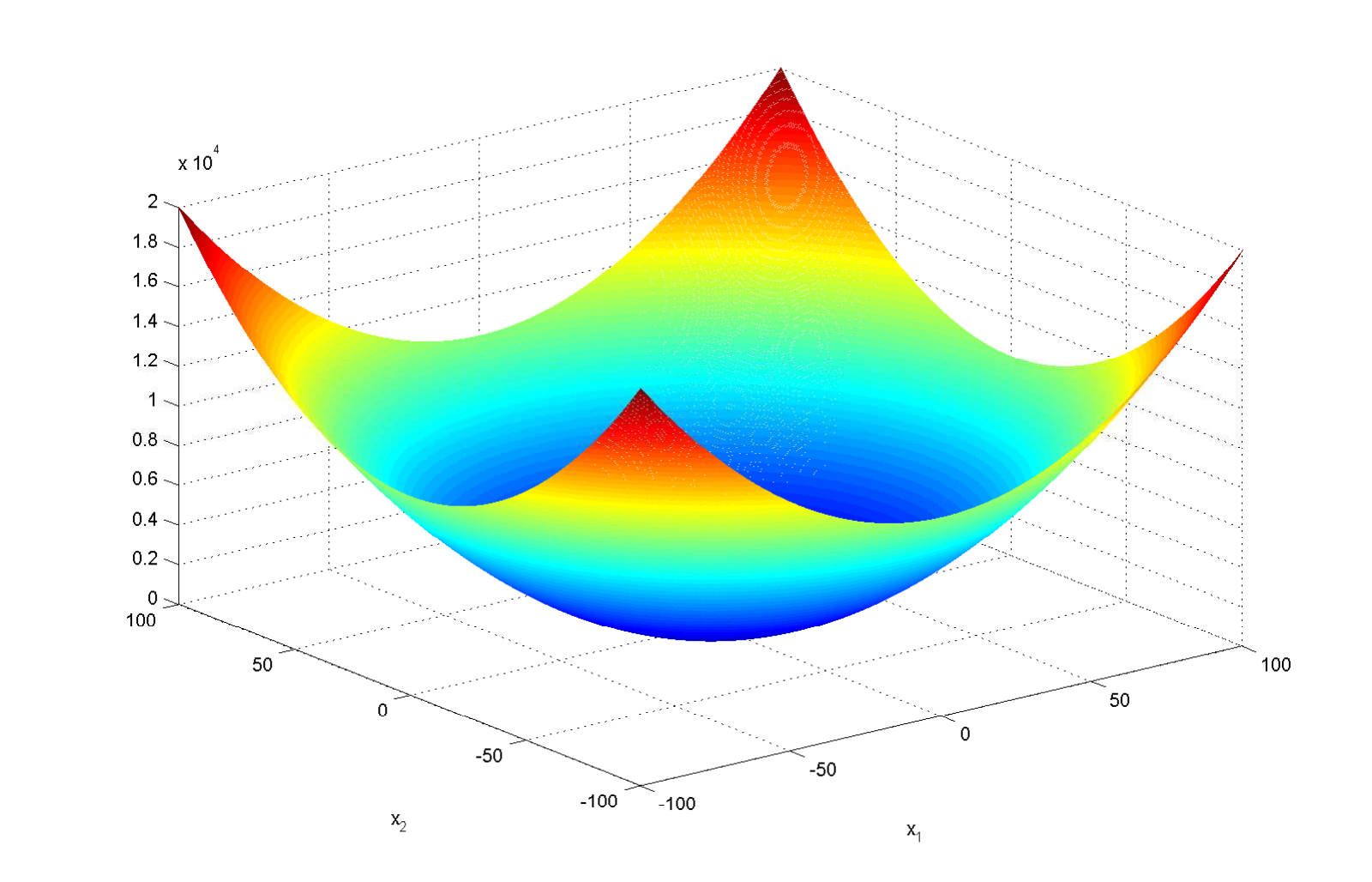}}
\subfloat[$f_2$]{\includegraphics[width=4cm,height=4cm]{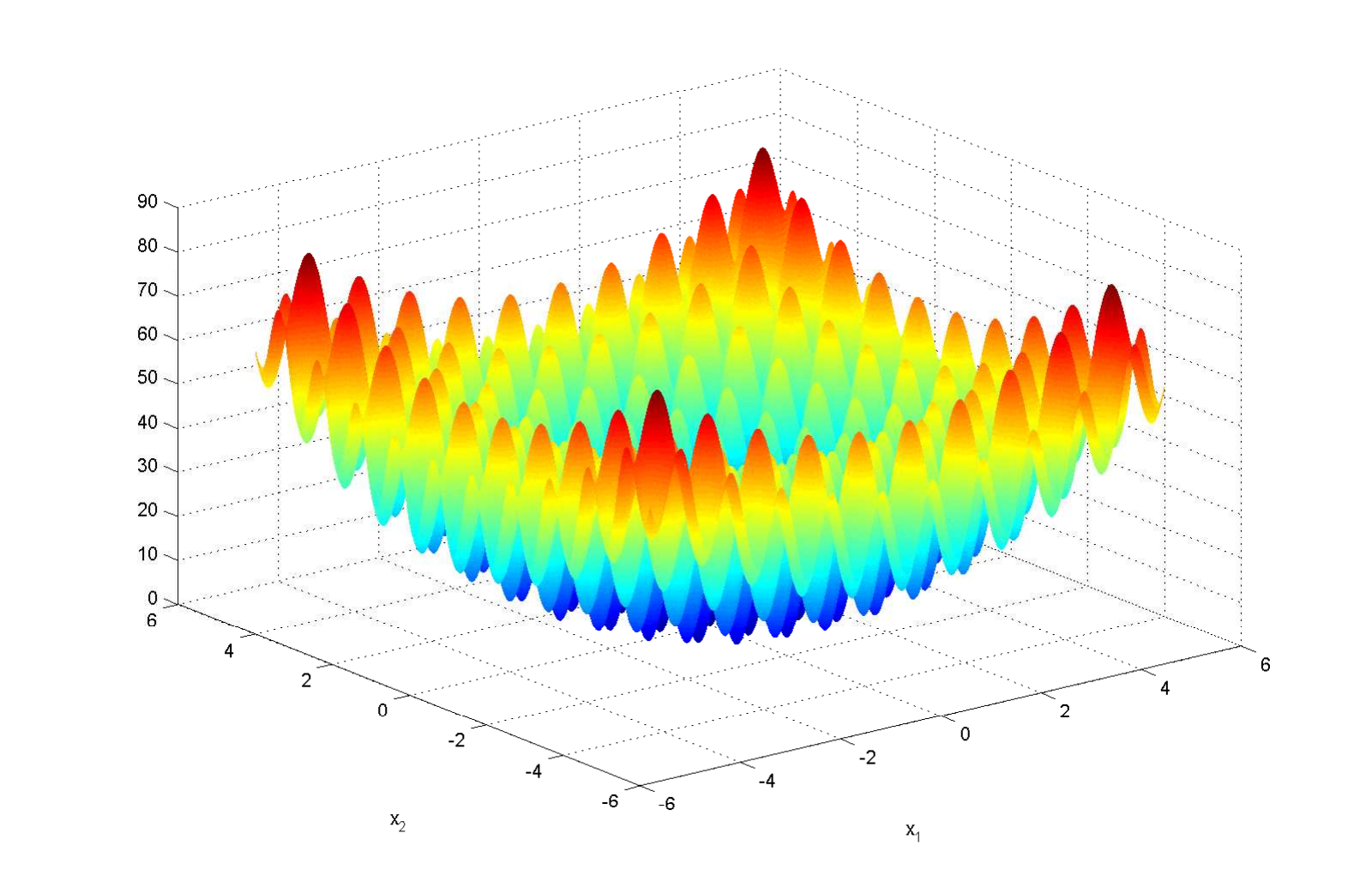}}
\subfloat[$f_3$]{\includegraphics[width=4cm,height=4cm]{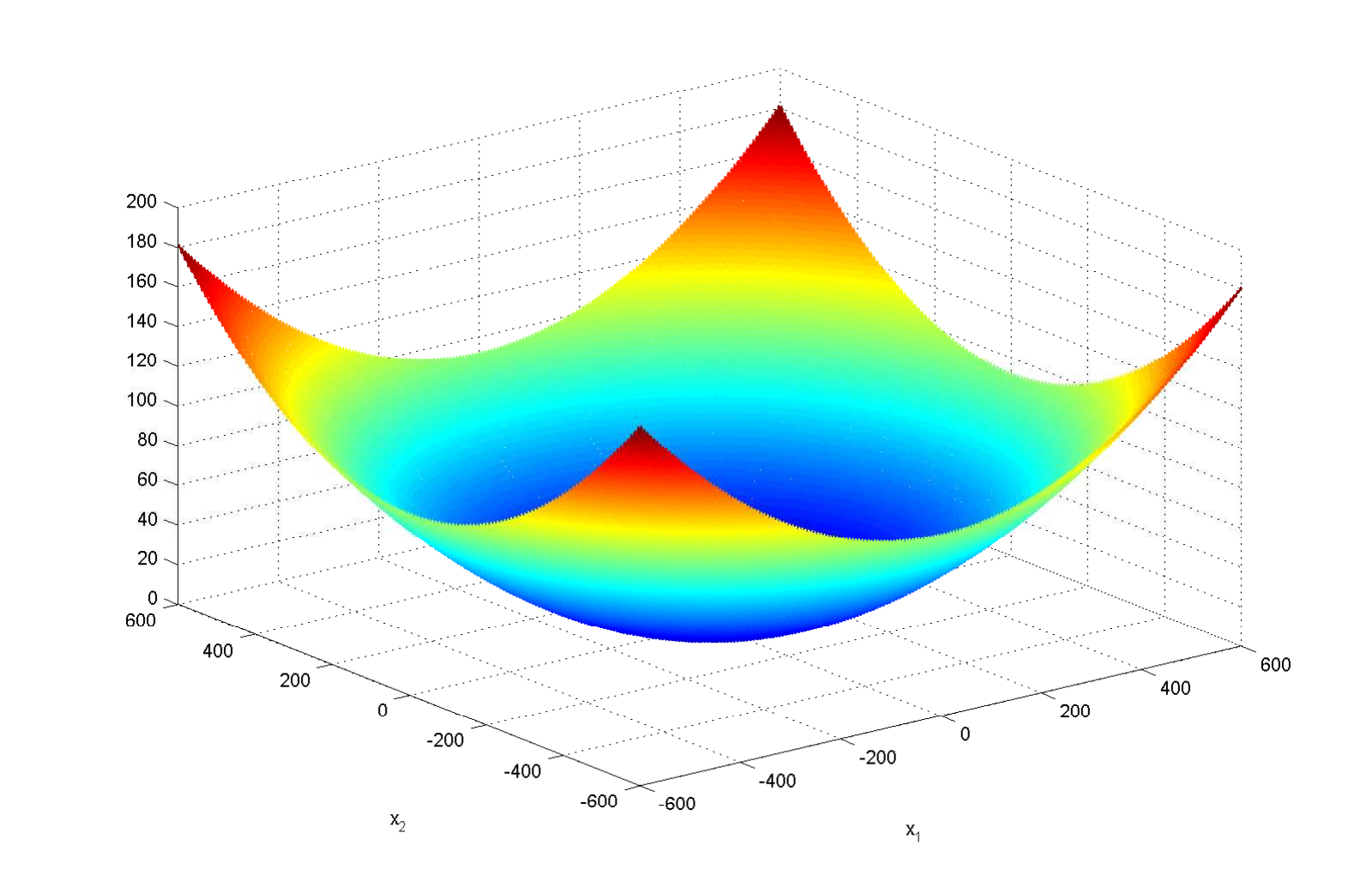}}\\
\subfloat[$f_4$]{\includegraphics[width=4cm,height=4cm]{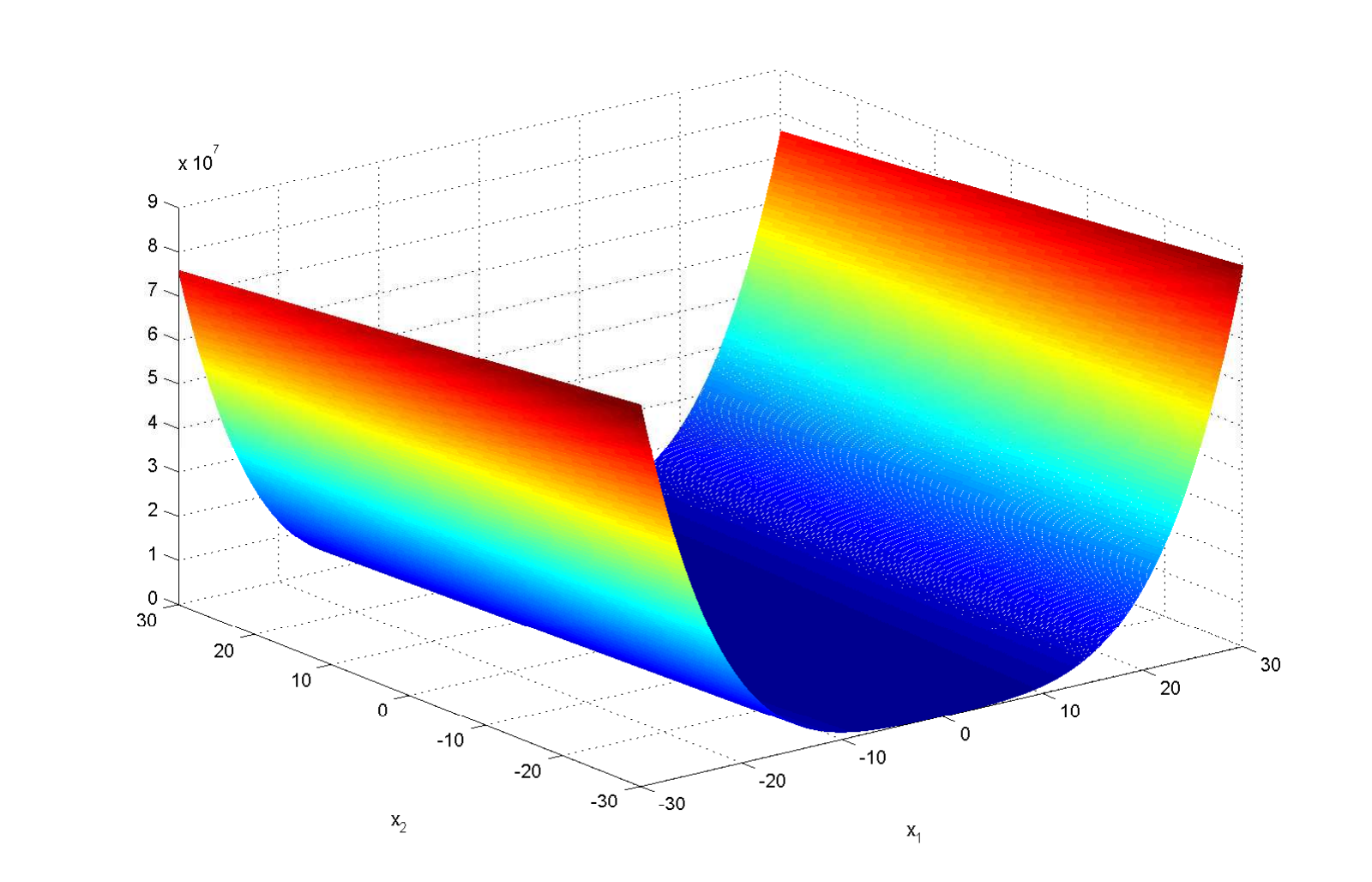}}
\subfloat[$f_5$]{\includegraphics[width=4cm,height=4cm]{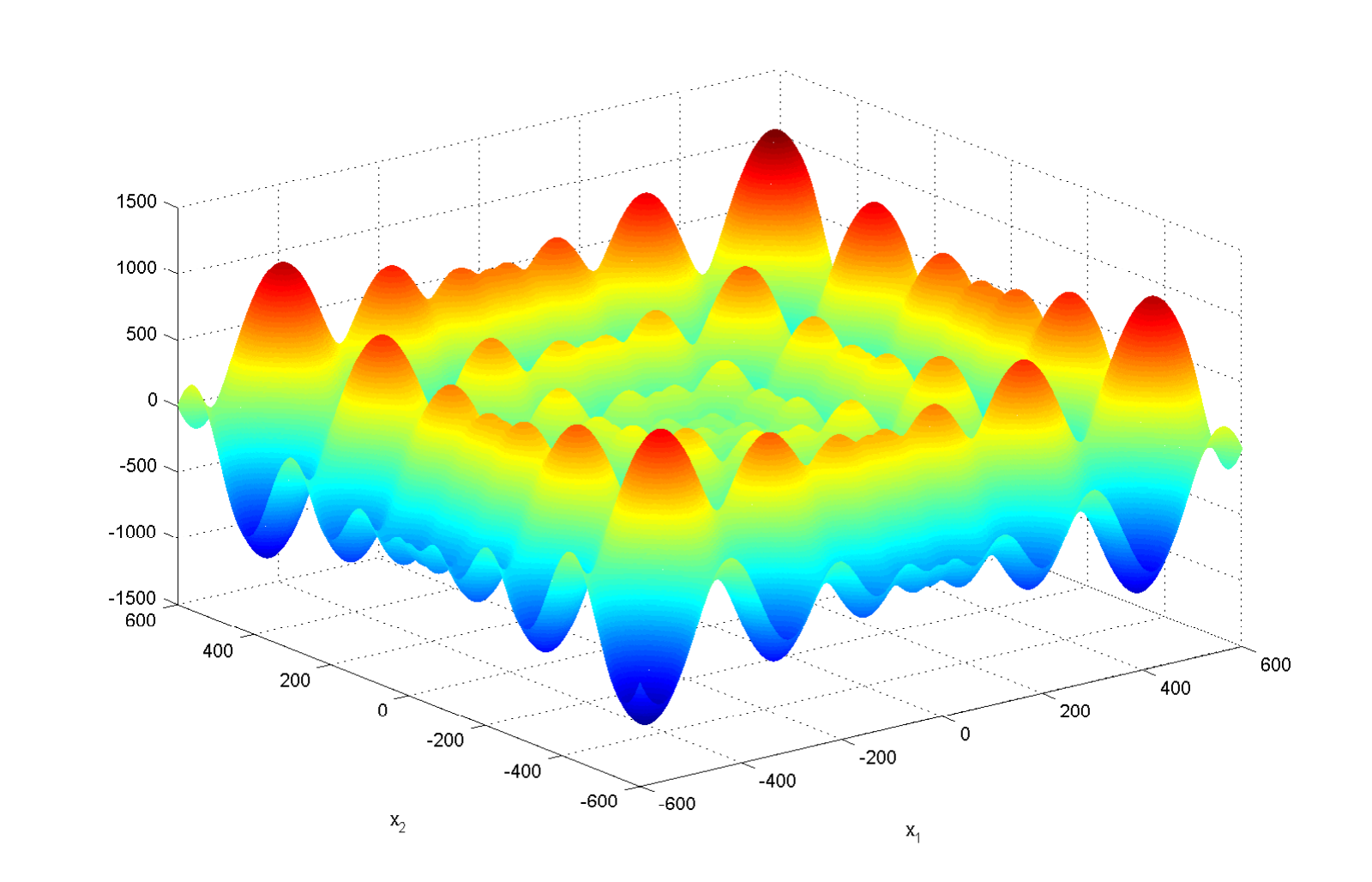}}
\subfloat[$f_6$]{\includegraphics[width=4cm,height=4cm]{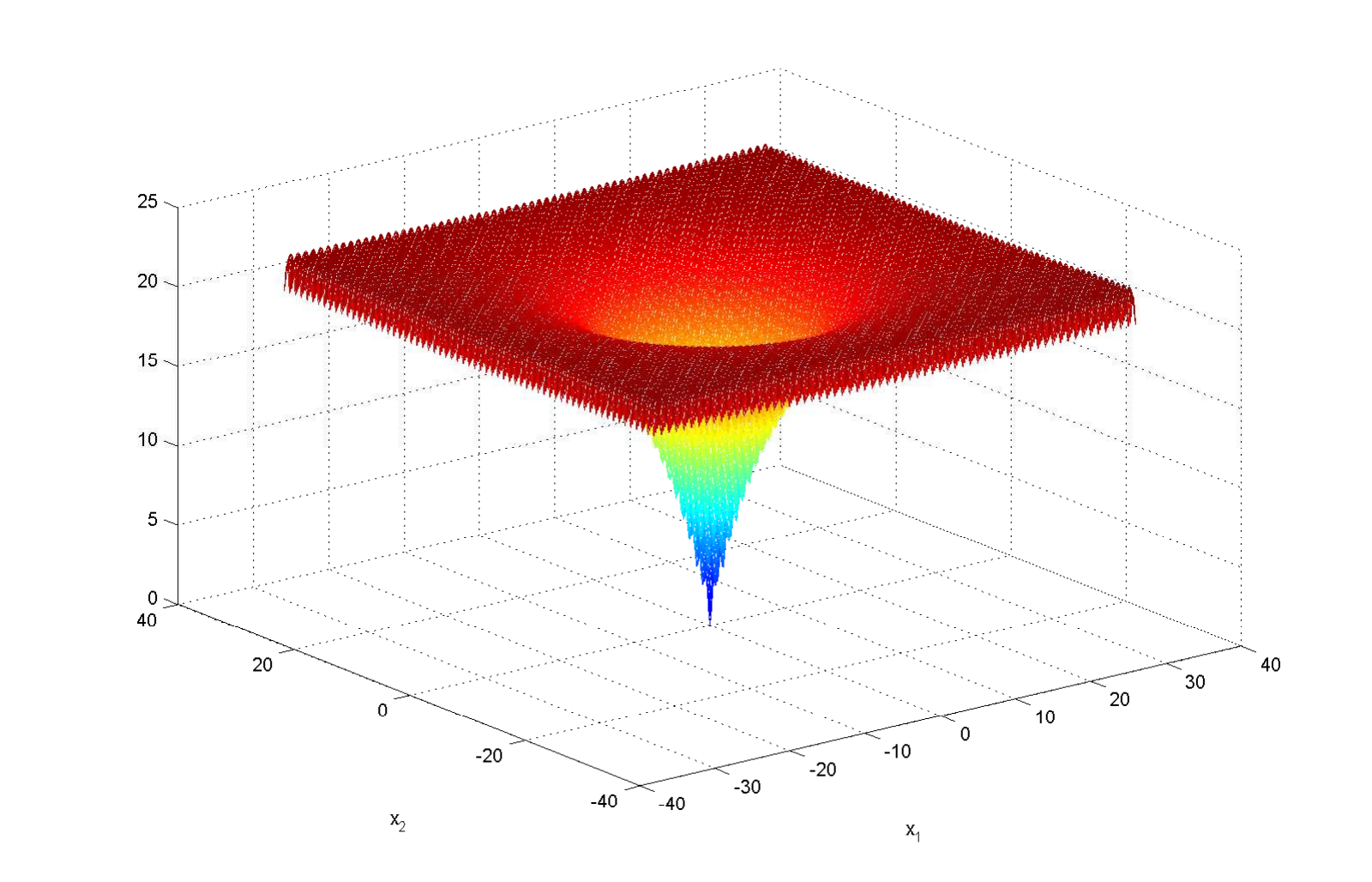}}\\
\subfloat[$f_7$]{\includegraphics[width=4cm,height=4cm]{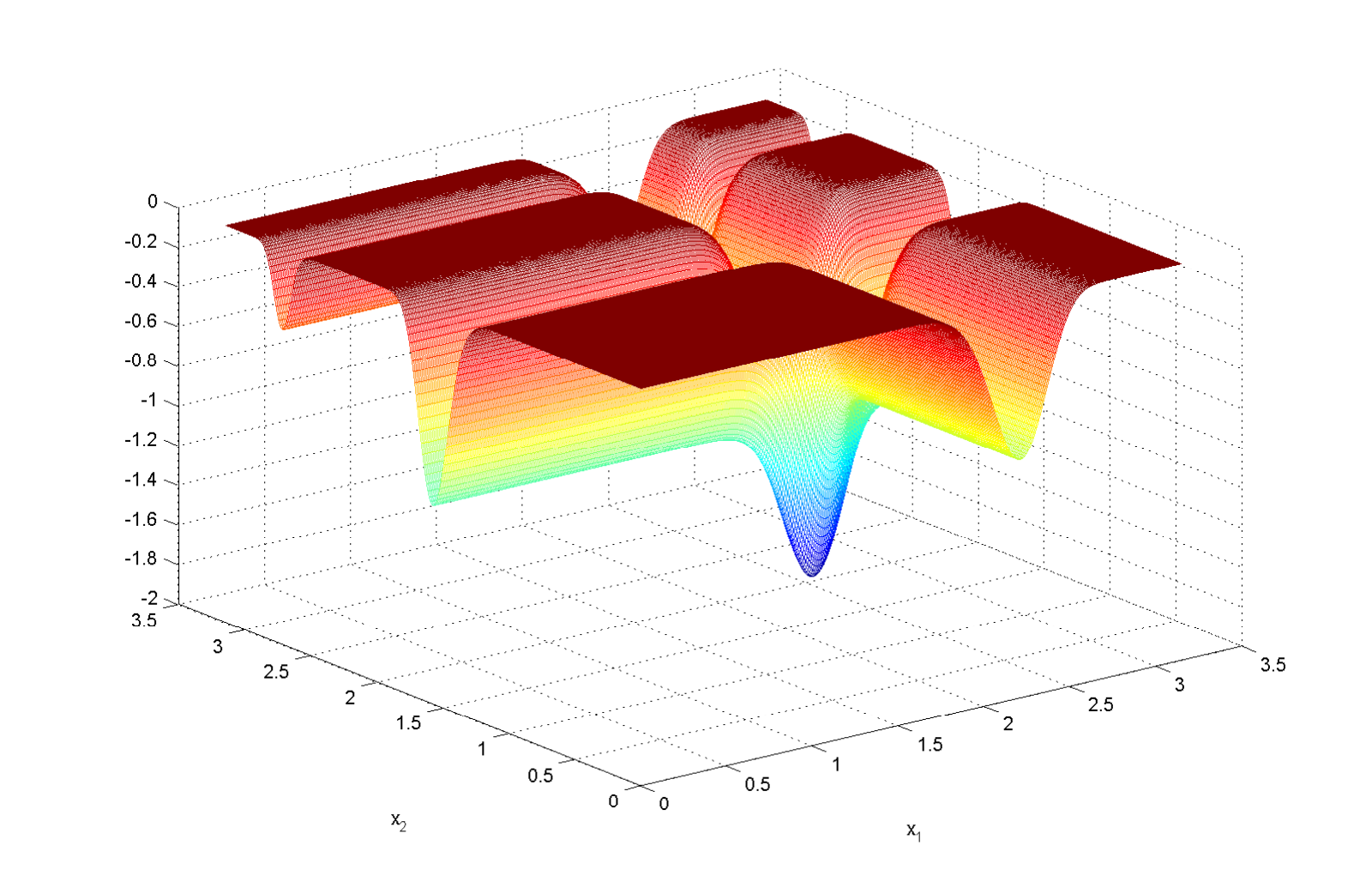}}
\subfloat[$f_8$]{\includegraphics[width=4cm,height=4cm]{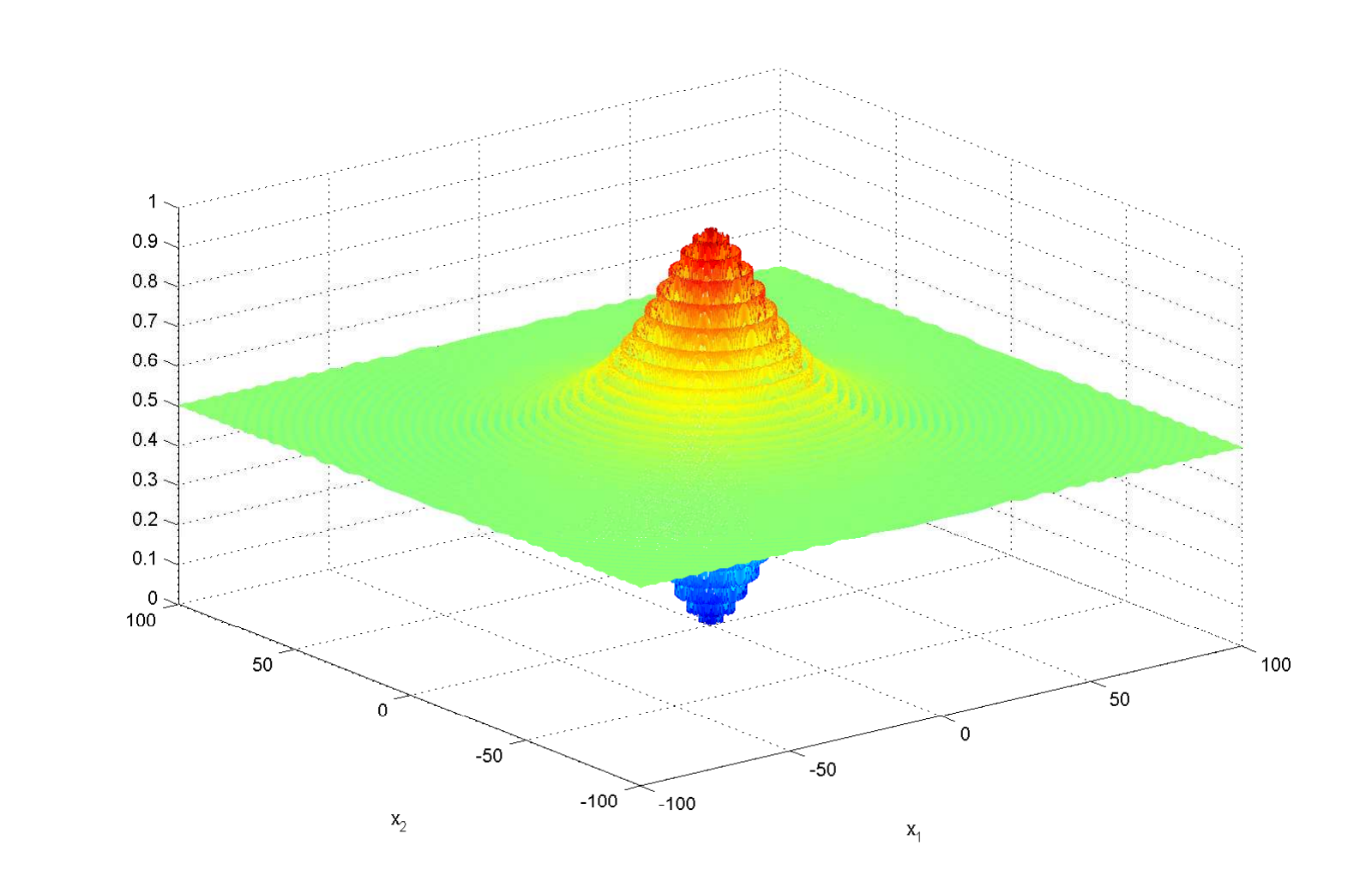}}
\subfloat[$f_9$]{\includegraphics[width=4cm,height=4cm]{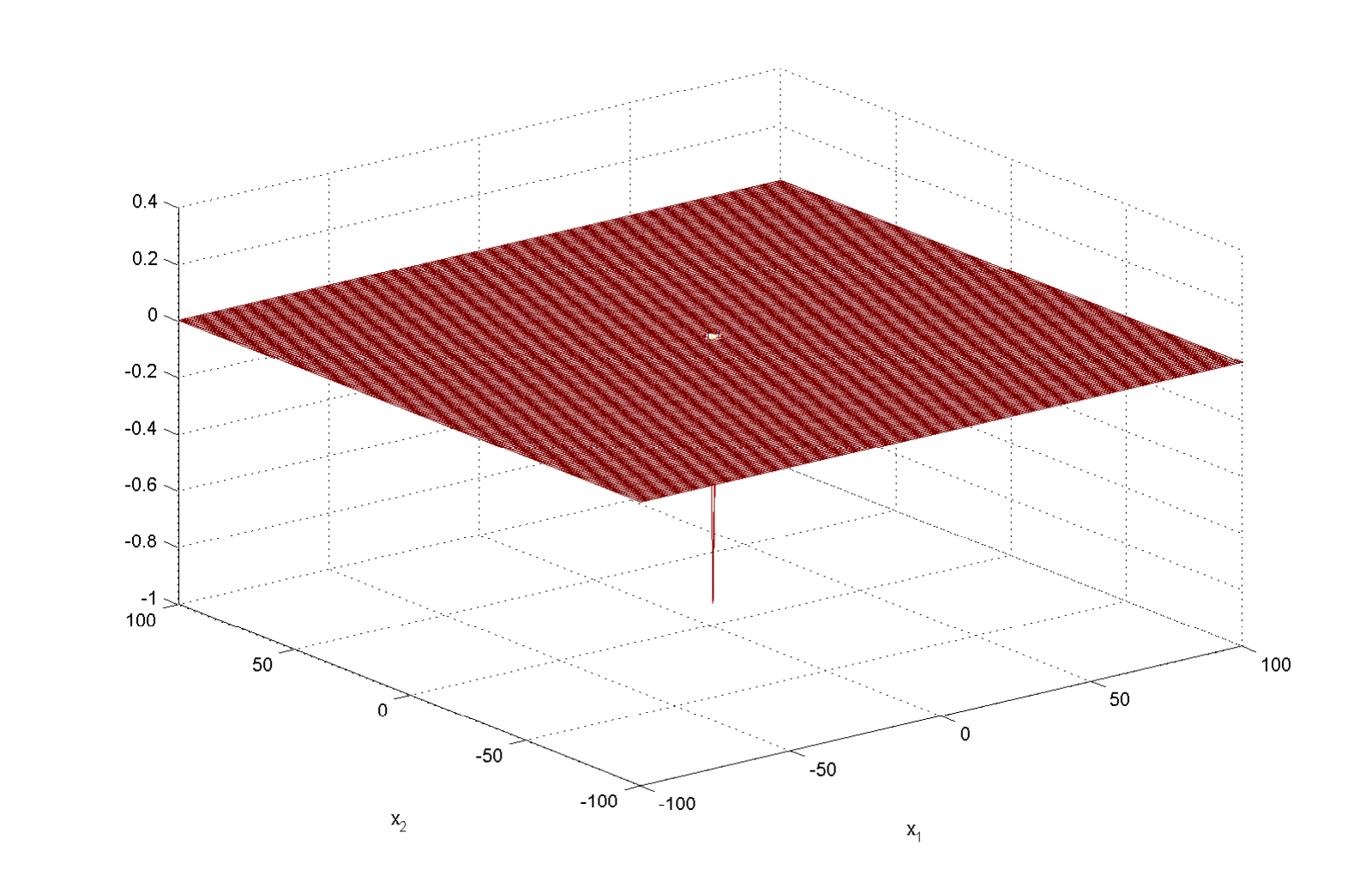}}\\
\subfloat[$f_{10}$]{\includegraphics[width=4cm,height=4cm]{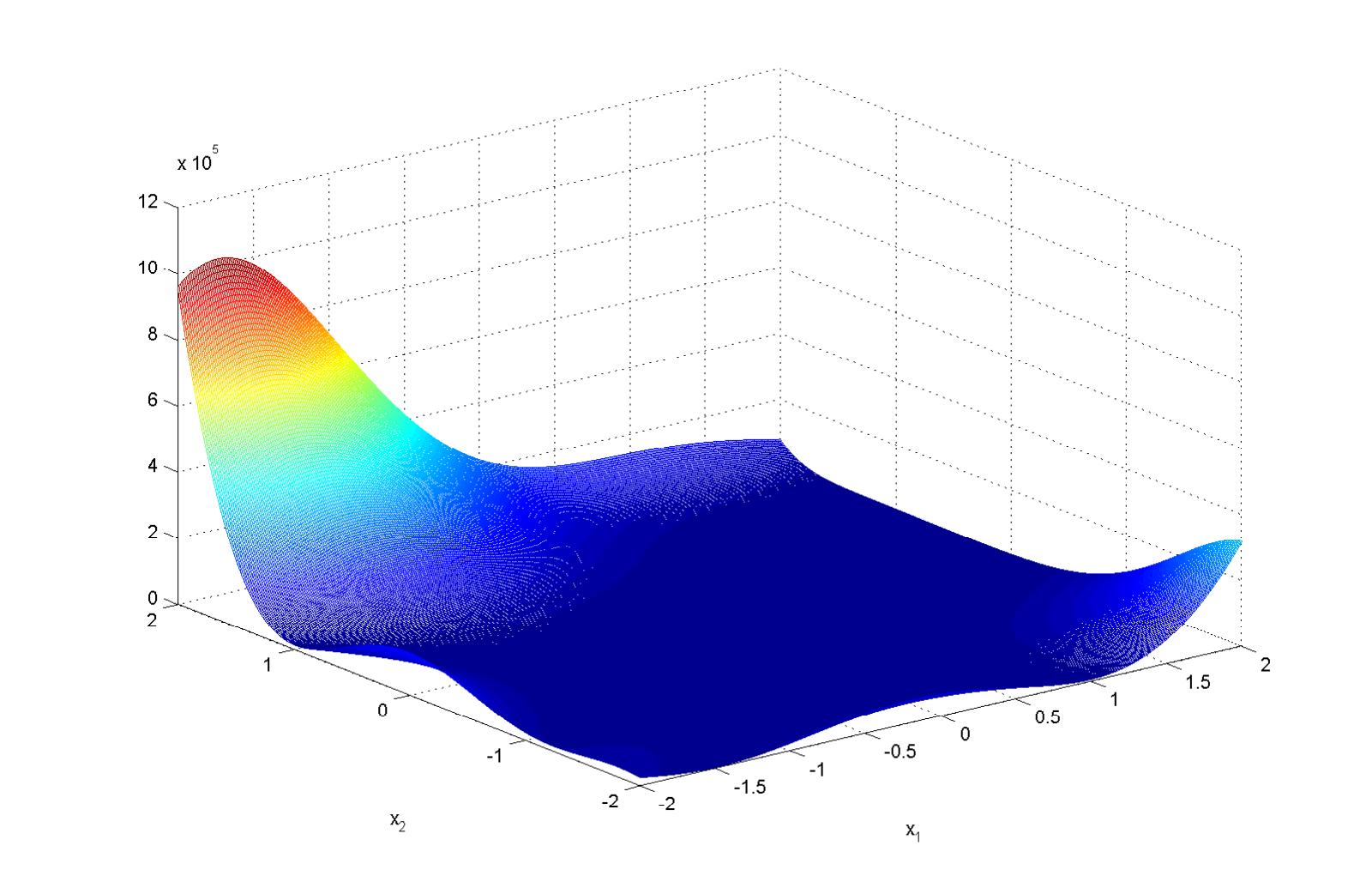}}
\caption{Plots of the two dimensional functions from $f_1$ to $f_{10}$}
\end{figure}
\subsection{Parameters setting}
All of the algorithms were run on MATLAB (Version R2010b) software platform. For simplicity and normalization, the experiment specifies all the control parameters of transformation operators in STAI and STAII starting at 1. Commonly, the variation of a parameter follows a linear, exponential or a logistic way. In this paper, the exponential way is accepted for its rapidity, of which the base is 2 in the experiment. In view of the operational precision of MATLAB in short format, the minimum $\alpha$ factor fixed at 1e-4 is enough for the situation.\\
\indent As for RCGA, we use the same parameter settings as in\cite{bi17}. Then, for CLPSO and SaDE, we use the MATLAB codes provided by the author in\cite{bi18,bi19} with minor revisions for this experiment.\\
\indent Programs were run independently for 30 trails, and for each trail, the population scale is 30, and the maximum iteration is 1000. The detailed parameters of STAI and STAII are shown in Table 2 and Table 3, respectively.
\begin{table}[h!]
\begin{center}
\caption{Parameters setting of STAI}
\footnotesize
\begin{tabular}{{cc}}
\hline
\toprule[1pt]
\textit{Parameter} & \textit{Value}\\
\hline
SE                     & 30 \\
$\alpha$               & 1 $\rightarrow$ 1e-4 \\
$\beta$                & 1  \\
$\gamma$               & 1  \\
$\delta$               & 1  \\
fc                     & 2  \\
\bottomrule[1pt]
\hline
\end{tabular}
\end{center}
\end{table}
\begin{table}[h!]
\begin{center}
\caption{Parameters setting of STAII}
\footnotesize
\begin{tabular}{{cc}}
\hline
\toprule[1pt]
\textit{Parameter} & \textit{Value}\\
\hline
SN                          & 30   \\
SE                          & 10   \\
CF                          & 50   \\
$\alpha$                    & 1 $\rightarrow$ 1e-4  \\
$\beta$                     & 1    \\
$\gamma$                    & 1    \\
$\delta$                    & 1  \\
fc                          & 2  \\
\bottomrule[1pt]
\hline
\end{tabular}
\end{center}
\end{table}
\subsection{Results and discussion}
For comparison, some common statistics are introduced. The \textit{best} means the minimum of the results, the \textit{worst} indicates the maximum of the results,
and then it follows the \textit{mean}, \textit{median} and \textit{st.dev.}(standard deviation). In some way, these statistics are able to evaluate
the search ability and solution accuracy, reliability and convergence as well as stability. To be more specific, the \textit{best} indicates the global search ability and solution accuracy, the \textit{worst} and the \textit{mean} signify the reliability and convergence, while the \textit{median} and \textit{st.dev.} correspond to the stability.\\
\indent Results for two dimensional functions optimization are listed in Table 4, while results for ten dimensional functions optimization can be found in Table 5. On the other hand, illustrations of the average fitness in 30 simulations are given in Figure 3 and Figure 4 for two dimensional and ten dimensional functions, respectively. The average fitness curve can visually depict the search ability and convergence performance. In the following paragraphs the analysis of the results for each functions will be discussed separately.
\begin{table}[!htbp]
\begin{center}
\caption{Comparisons among various algorithms on test functions(2D)}
\footnotesize
\begin{tabular}{{p{0.5cm}ccccccc}}
\hline
\toprule[1pt]
Fcn & Statistic&  RCGA & CLPSO & SaDE & STAI & STAII  \\
\hline
        & best	  & 1.5795e-028 &  2.7344e-091  & 2.9229e-196 & 0	        & 0           \\
        & median  & 4.9090e-024 &  2.6808e-087  & 4.2548e-189 & 0           & 0           \\
$f_{1}$ & mean    & 1.5458e-022 &  5.9580e-082  & 6.6729e-188 & 0	        & 0           \\
        & worst   & 2.5849e-021 &  1.7771e-080  & 9.6210e-187 & 0	        & 0           \\
        & st.dev. & 4.7939e-022 &  3.2440e-081  & 0           & 0	        & 0           \\
\hline
        & best	  & 0           &	0           & 0           & 0	        & 0           \\
        & median  & 0           &	0           & 0           & 0           & 0           \\
$f_{2}$ & mean    & 0           &	0	        & 0	          & 0	        & 0           \\
        & worst   & 0           &   0	        & 0	          & 0	        & 0           \\
        & st.dev. & 0           &   0	        & 0	          & 0	        & 0           \\
\hline
        & best	  & 0           &	0           & 0           & 0	        & 0           \\
        & median  & 0.0074      &	0           & 0           & 0           & 0           \\
$f_{3}$ & mean    & 0.0042      &	1.2460e-009	& 0	          & 0	        & 0           \\
        & worst   & 0.0074      &   3.7377e-008	& 0	          & 0	        & 0           \\
        & st.dev. & 0.0037      &   6.8241e-009	& 0	          & 0	        & 0           \\
\hline
        & best	  & 7.0832e-008 &  9.2890e-010  & 0           & 1.0092e-012	& 4.2592e-014 \\
        & median  & 0.0085      &  3.9984e-007  & 0           & 1.0900e-011 & 3.9400e-012 \\
$f_{4}$ & mean    & 1.0364      &  3.9260e-005	& 0	          & 1.2571e-011	& 4.4217e-012 \\
        & worst   & 26.2801     &  6.3652e-004	& 0	          & 4.7764e-011	& 1.4588e-011 \\
        & st.dev. & 4.7802      &  1.4088e-004  & 0	          & 1.0692e-011	& 3.9023e-012 \\
\hline
        & best	  & -837.9658   &	-837.9658   & -837.9658   & -837.9658	& -837.9658   \\
        & median  & -837.9658   &	-837.9658   & -837.9658   & -837.9658   & -837.9658   \\
$f_{5}$ & mean    & -822.1740   &	-837.9658	& -837.9658	  & -837.9658	& -837.9658   \\
        & worst   & -719.5274   &   -837.9658	& -837.9658	  & -837.9658	& -837.9658   \\
        & st.dev. & 40.9496     &   0	        & 0	          & 1.5939e-013	& 1.0970e-013 \\
\hline
        & best	  & 2.0428e-014 &	-8.8818e-016& -8.8818e-016& -8.8818e-016& -8.8818e-016\\
        & median  & 3.1681e-012 &	-8.8818e-016& -8.8818e-016& -8.8818e-016& -8.8818e-016\\
$f_{6}$ & mean    & 4.7516e-011 &	-8.8818e-016& -8.8818e-016& -8.8818e-016& -8.8818e-016\\
        & worst   & 9.6937e-010 &   -8.8818e-016& -8.8818e-016& -8.8818e-016& -8.8818e-016\\
        & st.dev. & 1.8319e-010 &   0	        & 0	          & 0	        & 0           \\
\hline
        & best	  & -1.8013     &	-1.8013     & -1.8013     & -1.8013	    &  -1.8013    \\
        & median  & -1.8013     &	-1.8013     & -1.8013     & -1.8013     &  -1.8013    \\
$f_{7}$ & mean    & -1.8013     &	-1.8013	    & -1.8013	  & -1.8013	    &  -1.8013    \\
        & worst   & -1.8013     &   -1.8013	    & -1.8013	  & -1.8013	    &  -1.8013    \\
        & st.dev. & 9.0336e-016 &   9.0336e-016	& 9.0336e-016 & 2.2063e-011	&  7.6618e-012\\
\hline
        & best	  & 0           &	0           & 0           & 0	        & 0           \\
        & median  & 0.0097      &	9.3299e-012 & 0           & 0           & 0           \\
$f_{8}$ & mean    & 0.0071      &	4.1836e-004	& 6.4773e-004 & 0	        & 0           \\
        & worst   & 0.0097      &   0.0097	    & 0.0097	  & 0	        & 0           \\
        & st.dev. & 0.0044      &   0.0018	    & 0.0025	  & 0	        & 0           \\
\hline
        & best	  & -1.0000     &	-1.0000     & -1.0000     & -1.0000	    & -1.0000     \\
        & median  & -1.0000     &	-1.0000     & -1.0000     & -1.0000     & -1.0000     \\
$f_{9}$ & mean    & -1.0000     &	-1.0000	    & -1.0000	  & -1.0000	    & -1.0000     \\
        & worst   & -1.0000     &   -1.0000	    & -1.0000	  & -1.0000	    & -1.0000     \\
        & st.dev. & 0           &   0	        & 0	          & 1.0124e-012	& 2.6511e-013 \\
\hline
        & best	  & 3.0000      &	3.0000      & 3.0000      & 3.0000	    & 3.0000      \\
        & median  & 3.0000      &	3.0000      & 3.0000      & 3.0000      & 3.0000      \\
$f_{10}$& mean    & 3.0000      &	3.0000	    & 3.0000	  & 3.0000	    & 3.0000      \\
        & worst   & 3.0000      &   3.0000	    & 3.0000	  & 3.0000	    & 3.0000      \\
        & st.dev. & 2.5135e-015 &   1.5317e-015	& 1.2669e-015 & 2.5697e-010	& 1.2191e-010 \\
\bottomrule[1pt]
\hline
\end{tabular}
\end{center}
\end{table}
\begin{figure}[!htbp]
\centering
\subfloat[$f_1$]{\includegraphics[width=4cm,height=4cm]{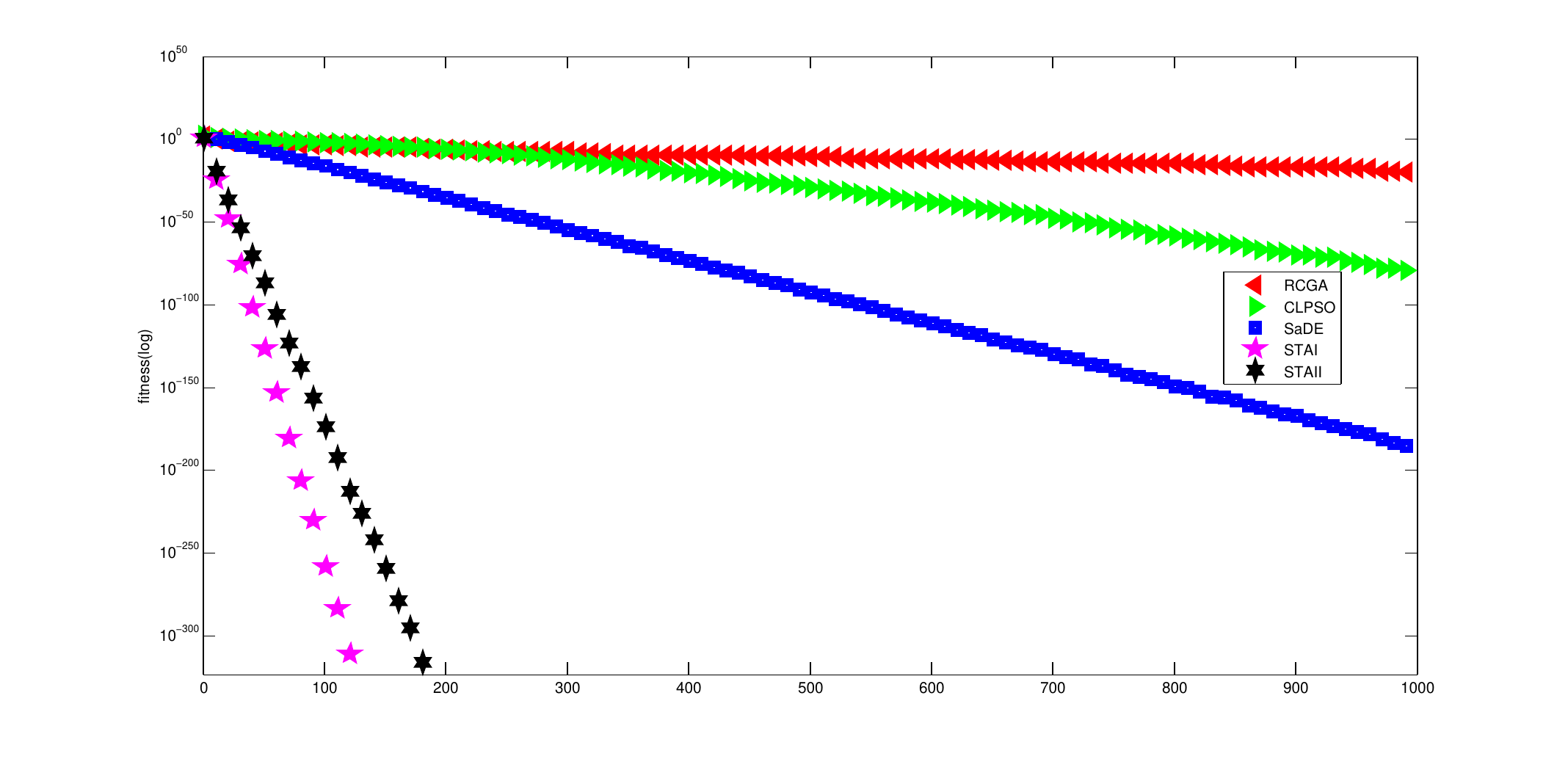}}
\subfloat[$f_2$]{\includegraphics[width=4cm,height=4cm]{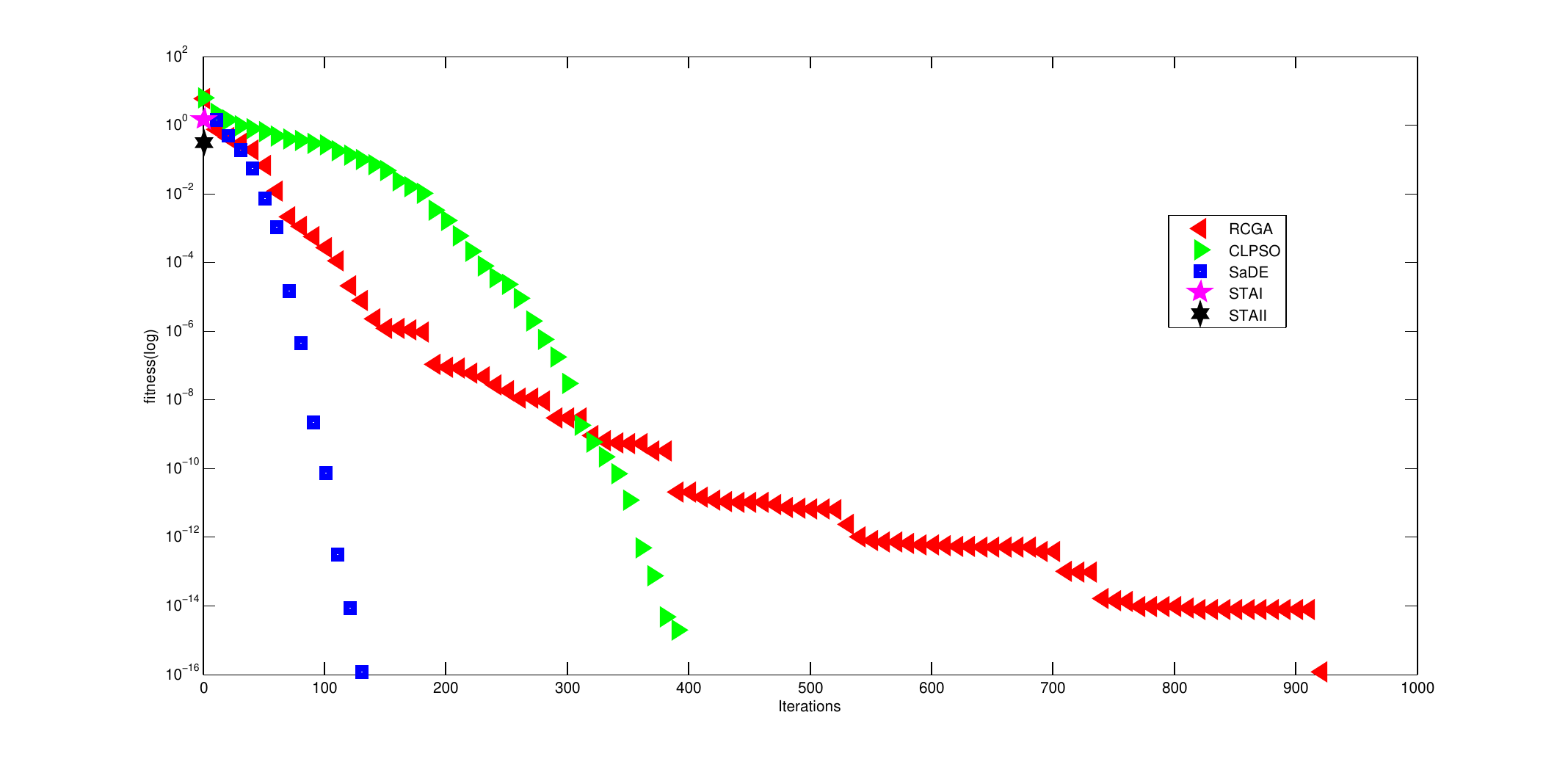}}
\subfloat[$f_3$]{\includegraphics[width=4cm,height=4cm]{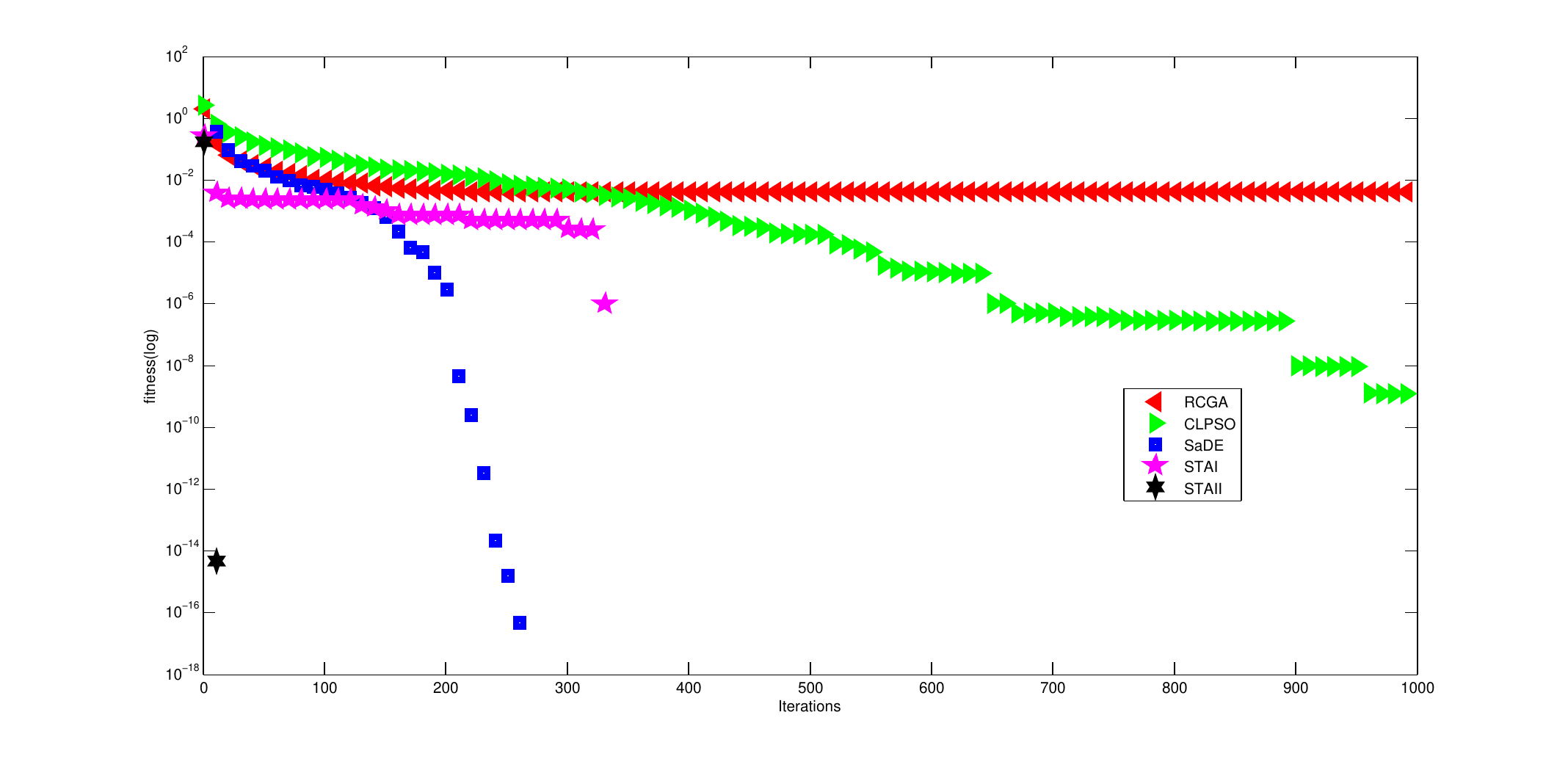}}\\
\subfloat[$f_4$]{\includegraphics[width=4cm,height=4cm]{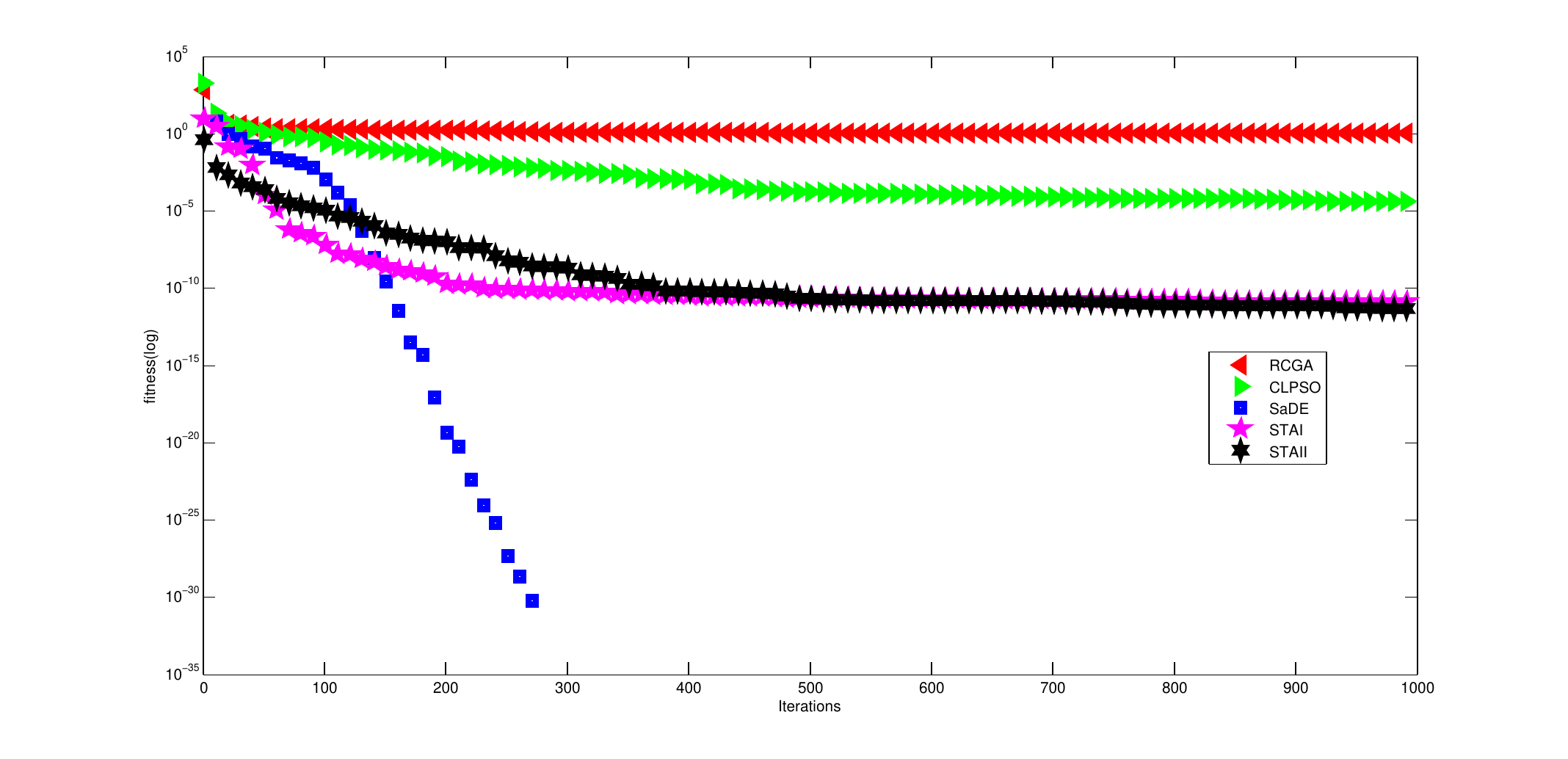}}
\subfloat[$f_5$]{\includegraphics[width=4cm,height=4cm]{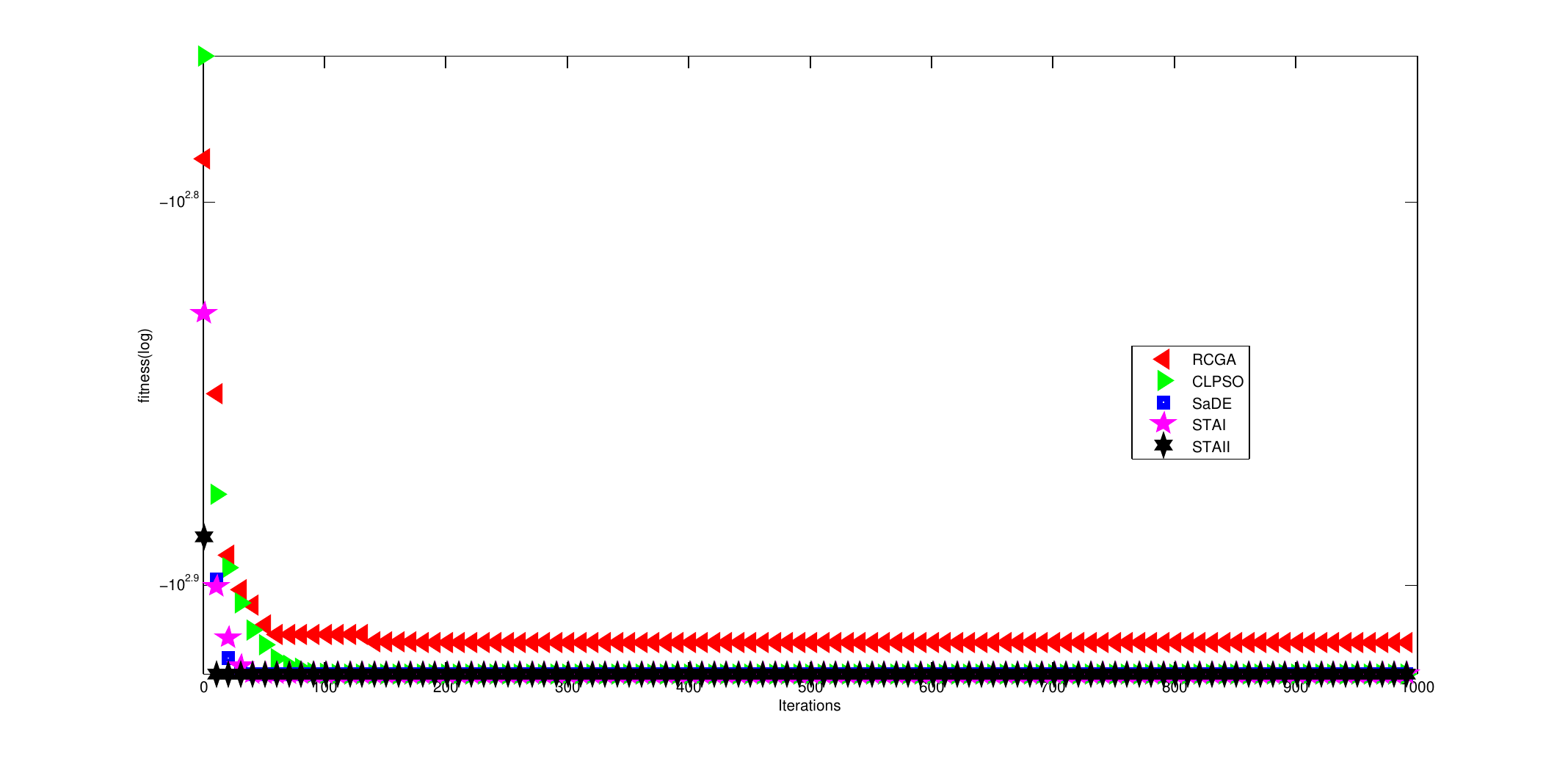}}
\subfloat[$f_6$]{\includegraphics[width=4cm,height=4cm]{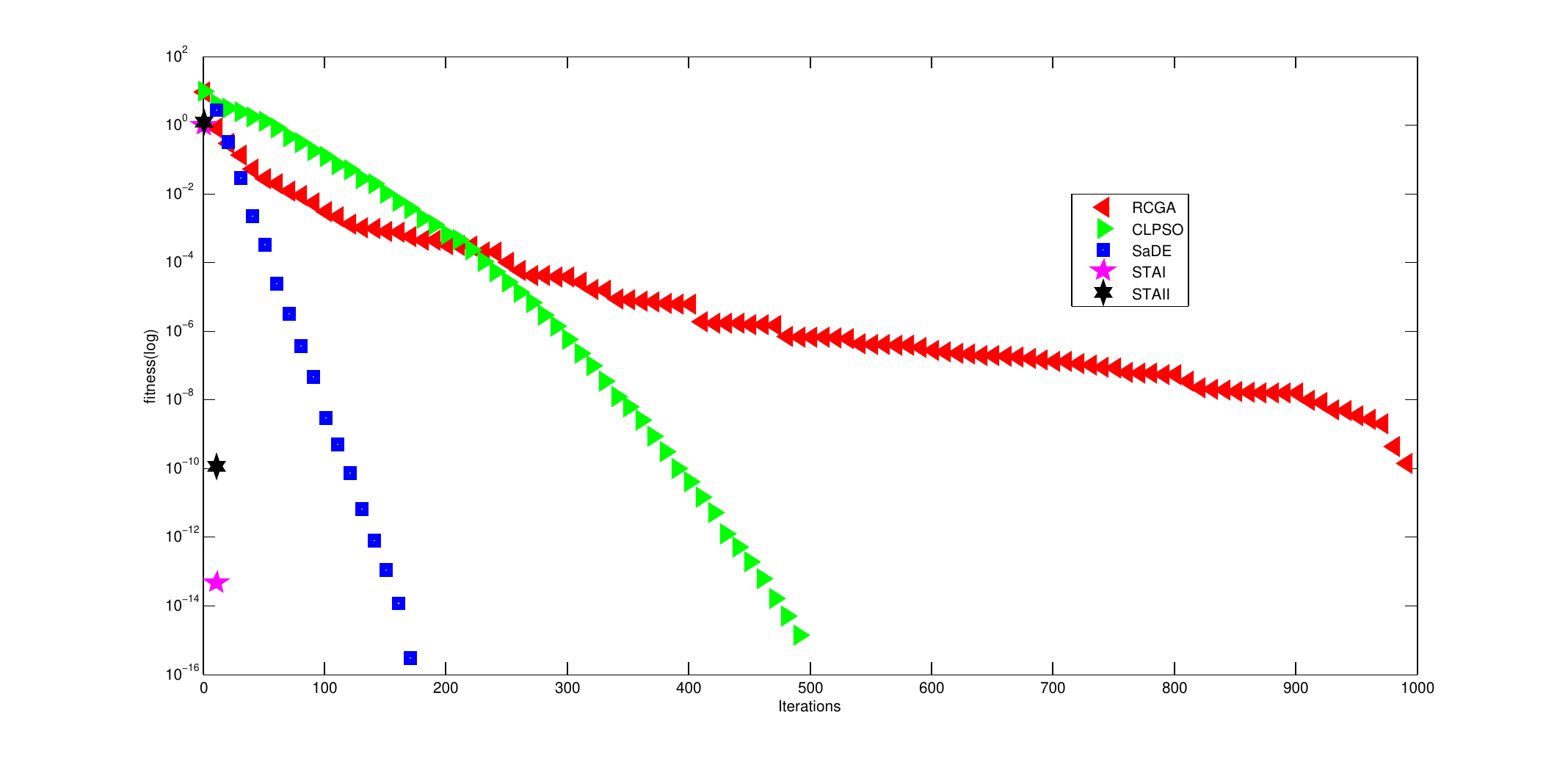}}\\
\subfloat[$f_7$]{\includegraphics[width=4cm,height=4cm]{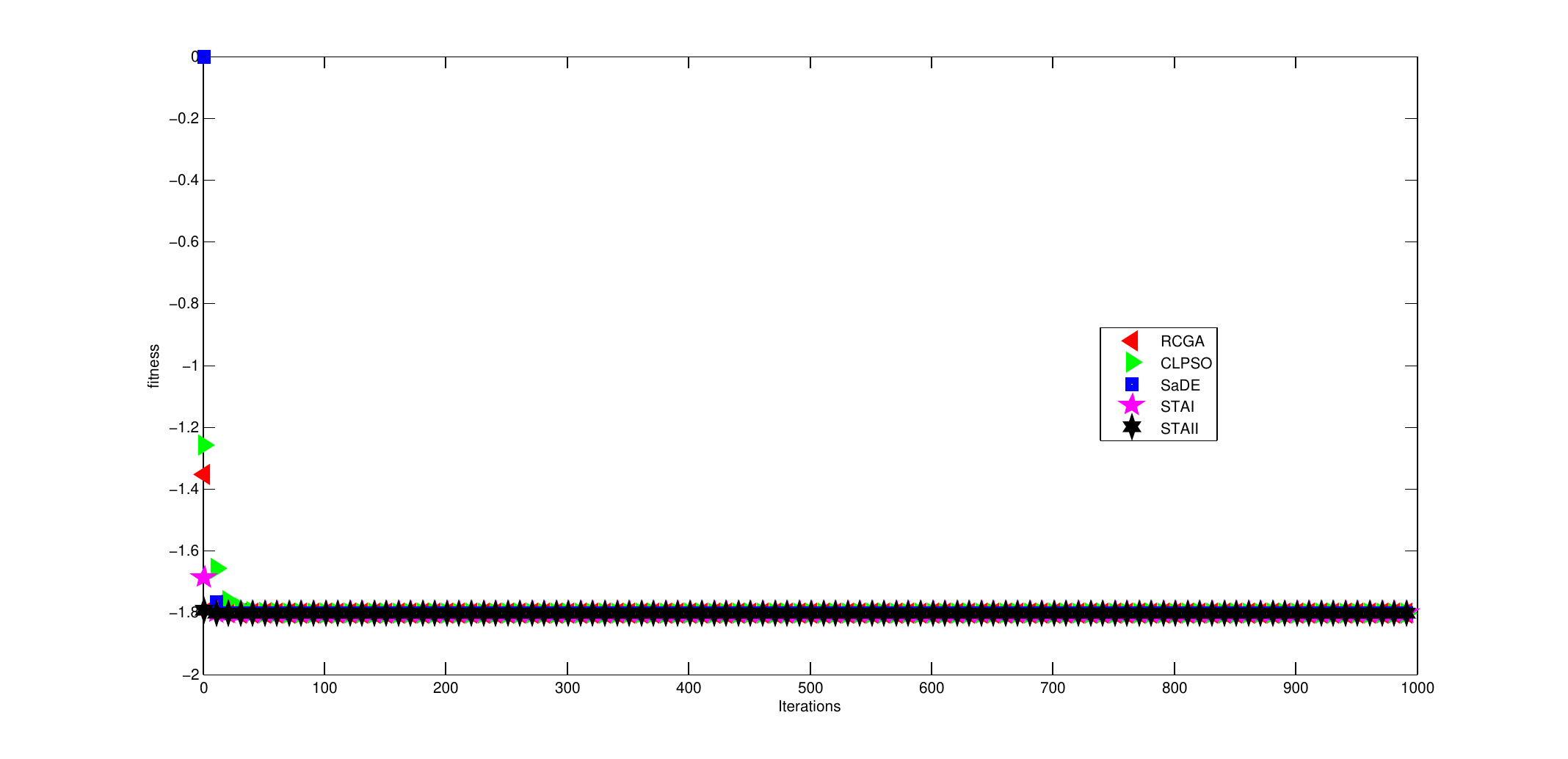}}
\subfloat[$f_8$]{\includegraphics[width=4cm,height=4cm]{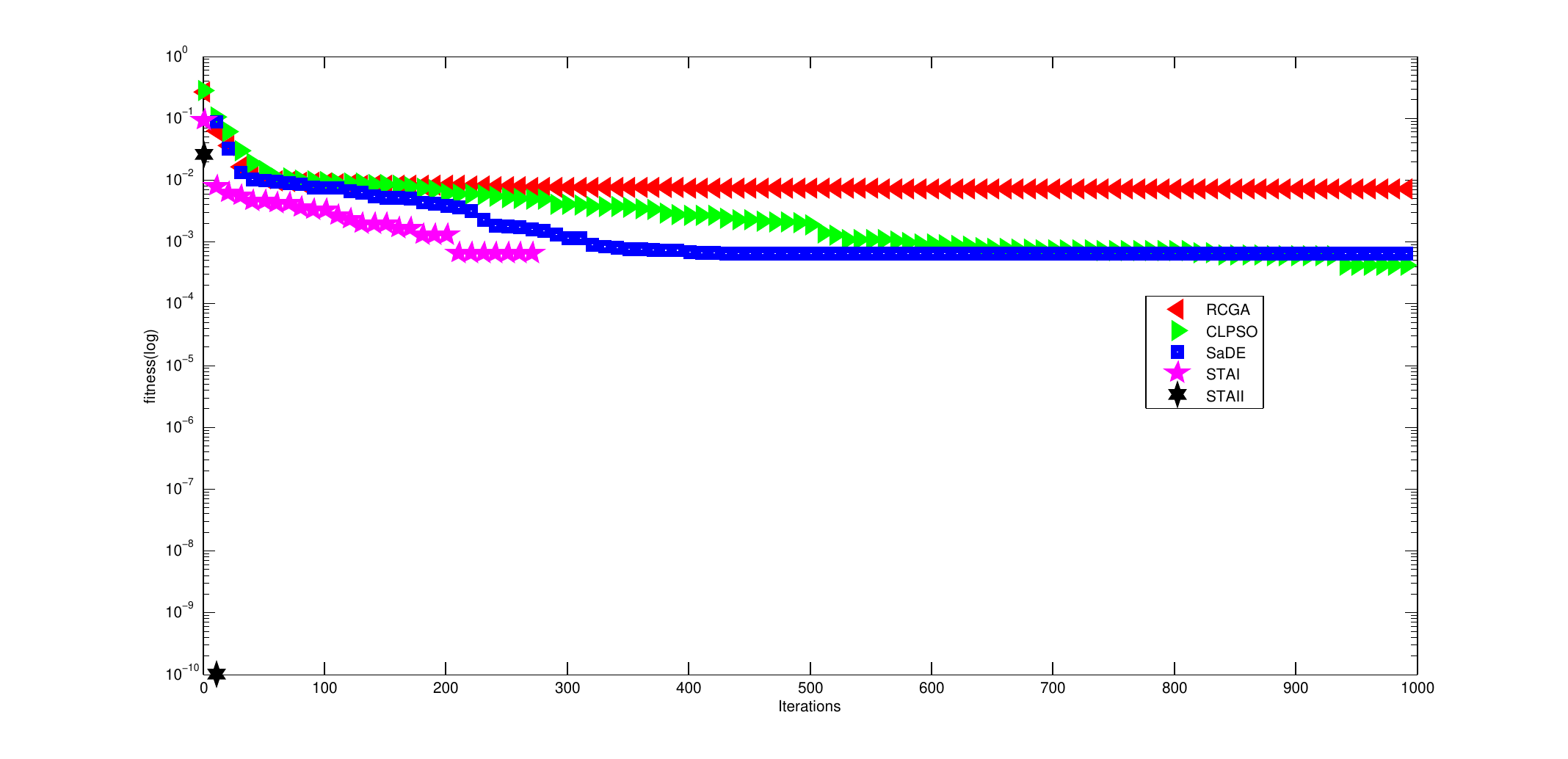}}
\subfloat[$f_9$]{\includegraphics[width=4cm,height=4cm]{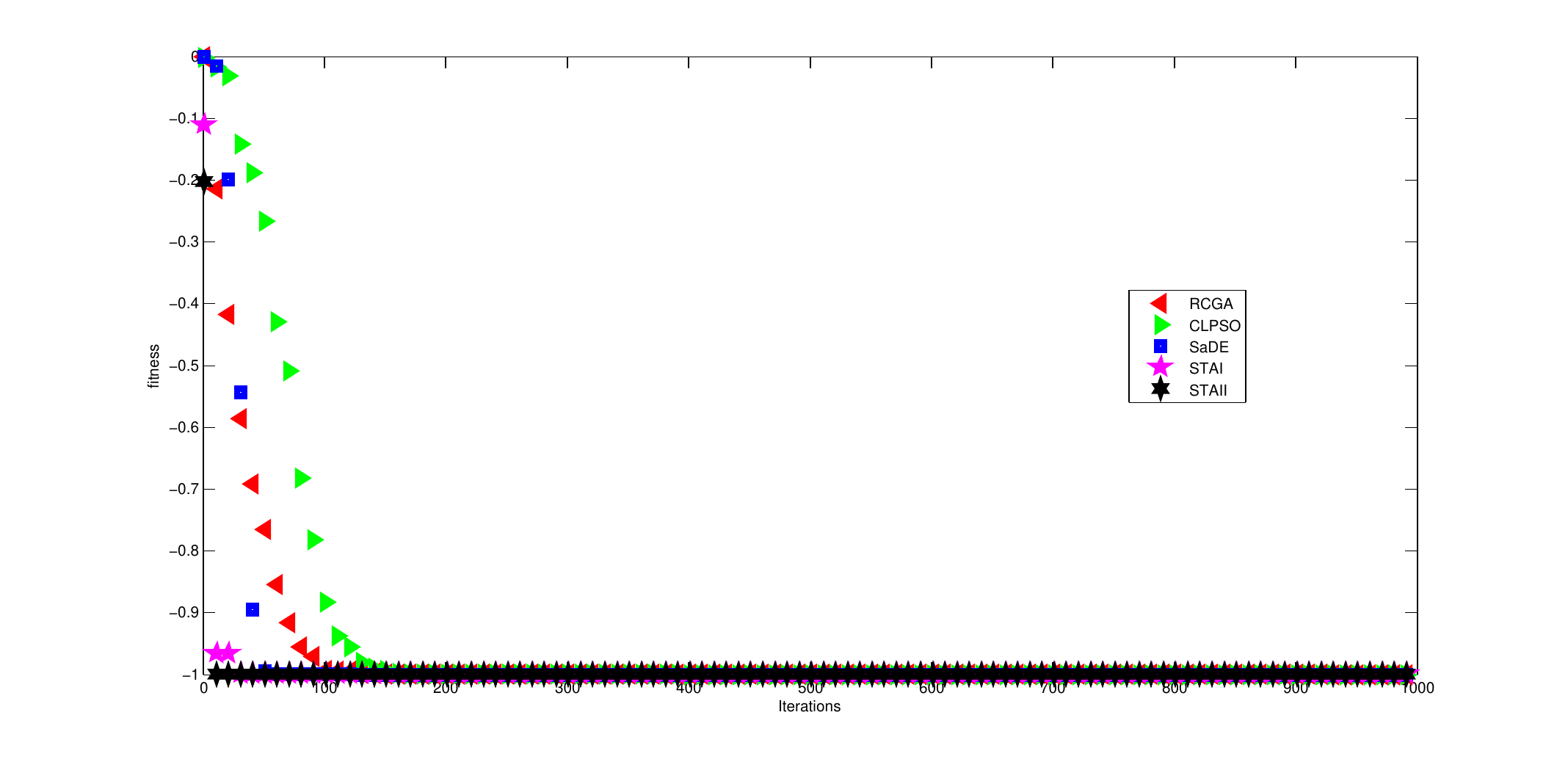}}\\
\subfloat[$f_{10}$]{\includegraphics[width=4cm,height=4cm]{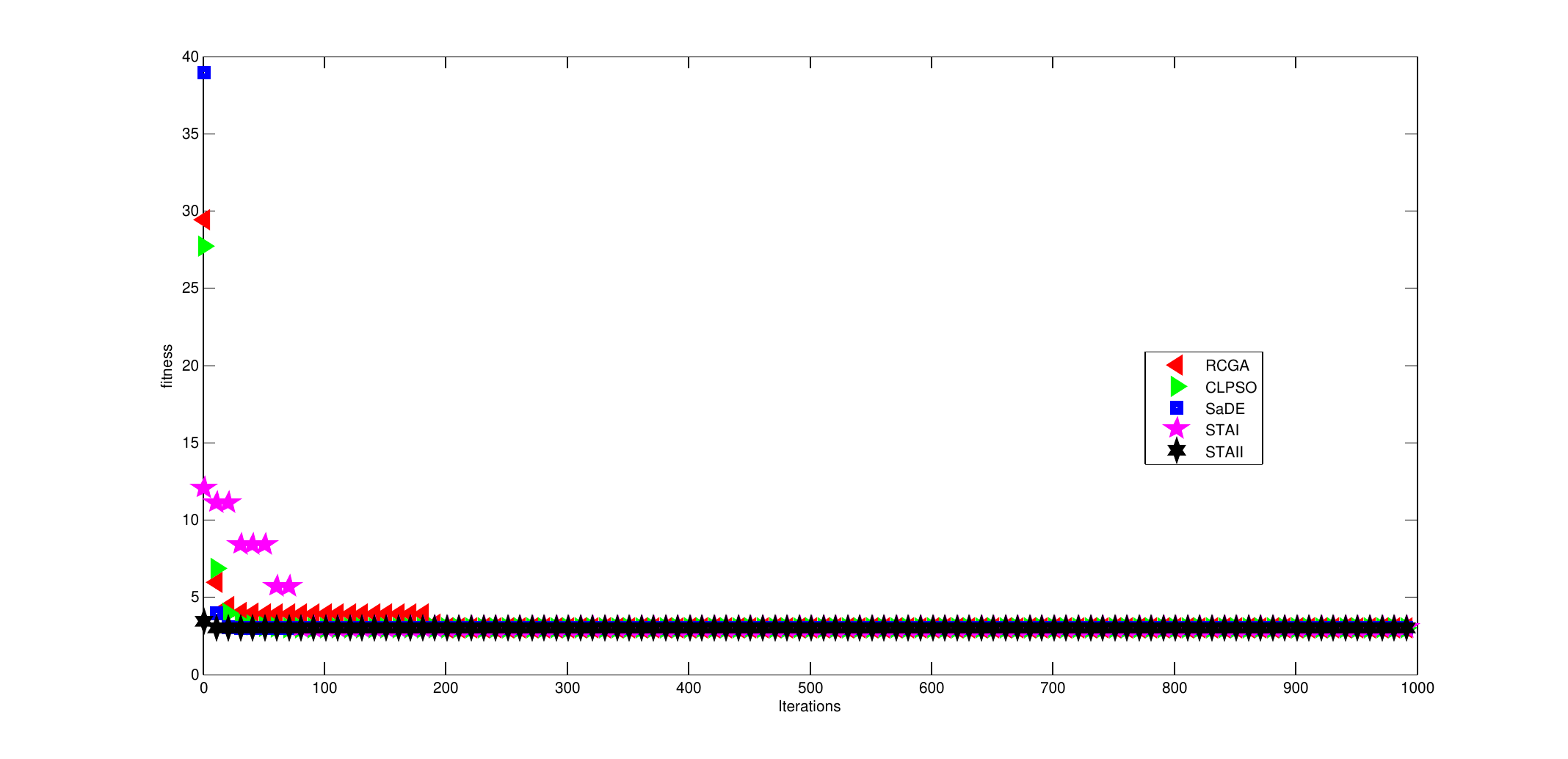}}
\caption{Average fitness of the two-dimensional functions from $f_1$ to $f_{10}$}
\end{figure}
\begin{table}[!htbp]
\begin{center}
\caption{Comparisons among various algorithms on test functions(10D)}
\footnotesize
\begin{tabular}{{p{0.5cm}ccccccc}}
\hline
\toprule[1pt]
Fcn & Statistic&  RCGA & CLPSO & SaDE & STAI & STAII  \\
\hline
        & best	  & 4.0118e-012  &	1.5229e-012  & 1.0686e-053 & 0	        & 0           \\
        & median  & 4.1019e-011  &	5.3464e-011  & 1.6873e-051 & 0          & 0           \\
$f_{1}$ & mean    & 8.4302e-011  &	6.0282e-011  & 3.5549e-051 & 0	        & 0           \\
        & worst   & 5.6332e-010  &  1.6279e-010  & 2.7214e-050 & 0	        & 0           \\
        & st.dev. & 1.1766e-010  &  4.9198e-011  & 6.1159e-051 & 0	        & 0           \\
\hline
        & best	  & 1.0814e-011  &	1.9806e-006  & 0           & 0	        & 0           \\
        & median  & 2.9849       &	1.0401e-005  & 0           & 0          & 0           \\
$f_{2}$ & mean    & 2.6864       &	2.7769e-005	 & 0.0332	   & 0	        & 0           \\
        & worst   & 5.9697       &  2.5009e-004	 & 0.9950	   & 0	        & 0           \\
        & st.dev. & 1.5491       &  4.7363e-005	 & 0.1817	   & 0	        & 0           \\
\hline
        & best	  & 3.2914e-010  &	2.6061e-005  & 0           & 0	        & 0           \\
        & median  & 0.0492       &	0.0019       & 0           & 0          & 0           \\
$f_{3}$ & mean    & 0.0582       &	0.0038	     & 5.7529e-004 & 0.0166	    & 0           \\
        & worst   & 0.1699       &  0.0166	     & 0.0099	   & 0.0738	    & 0           \\
        & st.dev. & 0.0439       &  0.0045	     & 0.0022	   & 0.0260	    & 0           \\
\hline
        & best	  & 0.0662       &	0.6841       & 8.3515e-012 & 2.6607e-005& 7.3949e-005 \\
        & median  & 7.2483       &	4.0404       & 3.4637e-004 & 1.7823     & 0.2809      \\
$f_{4}$ & mean    & 7.0110       &	4.3775	     & 0.3415	   & 2.3266	    & 0.4095      \\
        & worst   & 9.2754       &  12.3253	     & 3.9866	   & 21.8603	& 1.5228      \\
        & st.dev. & 1.4466       &  2.8042       & 1.0098	   & 3.7249	    & 0.4124      \\
\hline
        & best	  & -3.9530e+003 &	-4.1898e+003 & -4.1898e+003&  -4.1898e+003& -4.1898e+003\\
        & median  & -3.8345e+003 &	-4.1898e+003 & -4.1898e+003&  -4.1898e+003& -4.1898e+003\\
$f_{5}$ & mean    & -3.7832e+003 &	-4.1898e+003 & -4.1898e+003&  -4.1898e+003& -4.1898e+003\\
        & worst   & -3.3608e+003 &  -4.1898e+003 & -4.1898e+003&  -4.1898e+003& -4.1898e+003\\
        & st.dev. & 151.3665     &  5.5006e-008	 & 2.7751e-012 &  1.1734e-011 & 1.9256e-012 \\
\hline
        & best	  & 6.6339e-007  &	1.5533e-006  & -8.8818e-016& -8.8818e-016 & -8.8818e-016\\
        & median  & 2.8352e-006  &	4.0369e-006  & 2.6645e-015 & -8.8818e-016 & -8.8818e-016\\
$f_{6}$ & mean    & 3.1524e-006  &	4.8197e-006  & 2.3093e-015 & 2.9606e-016  & -8.8818e-016\\
        & worst   & 9.6980e-006  &  1.3052e-005  & 2.6645e-015 & 2.6645e-015  & -8.8818e-016\\
        & st.dev. & 2.1010e-006  &  3.2984e-006	 & 1.0840e-015 & 1.7034e-015  & 0           \\
\hline
        & best	  & -9.6154      &	-9.6601      & -9.6602     & -9.6602	  &  -9.6602    \\
        & median  & -9.2604      &	-9.6598      & -9.6602     & -9.6602      &  -9.6602    \\
$f_{7}$ & mean    & -9.2425      &	-9.6588	     & -9.6513	   & -9.1797	  &  -9.6602    \\
        & worst   & -8.7143      &  -9.6549	     & -9.6135	   & -7.6602	  &  -9.6602    \\
        & st.dev. & 0.2265       &  0.0018	     & 0.0173      & 0.6044       &  1.9138e-009\\
\bottomrule[1pt]
\hline
\end{tabular}
\end{center}
\end{table}
\begin{figure}[!htbp]
\centering
\subfloat[$f_1$]{\includegraphics[width=4cm,height=4cm]{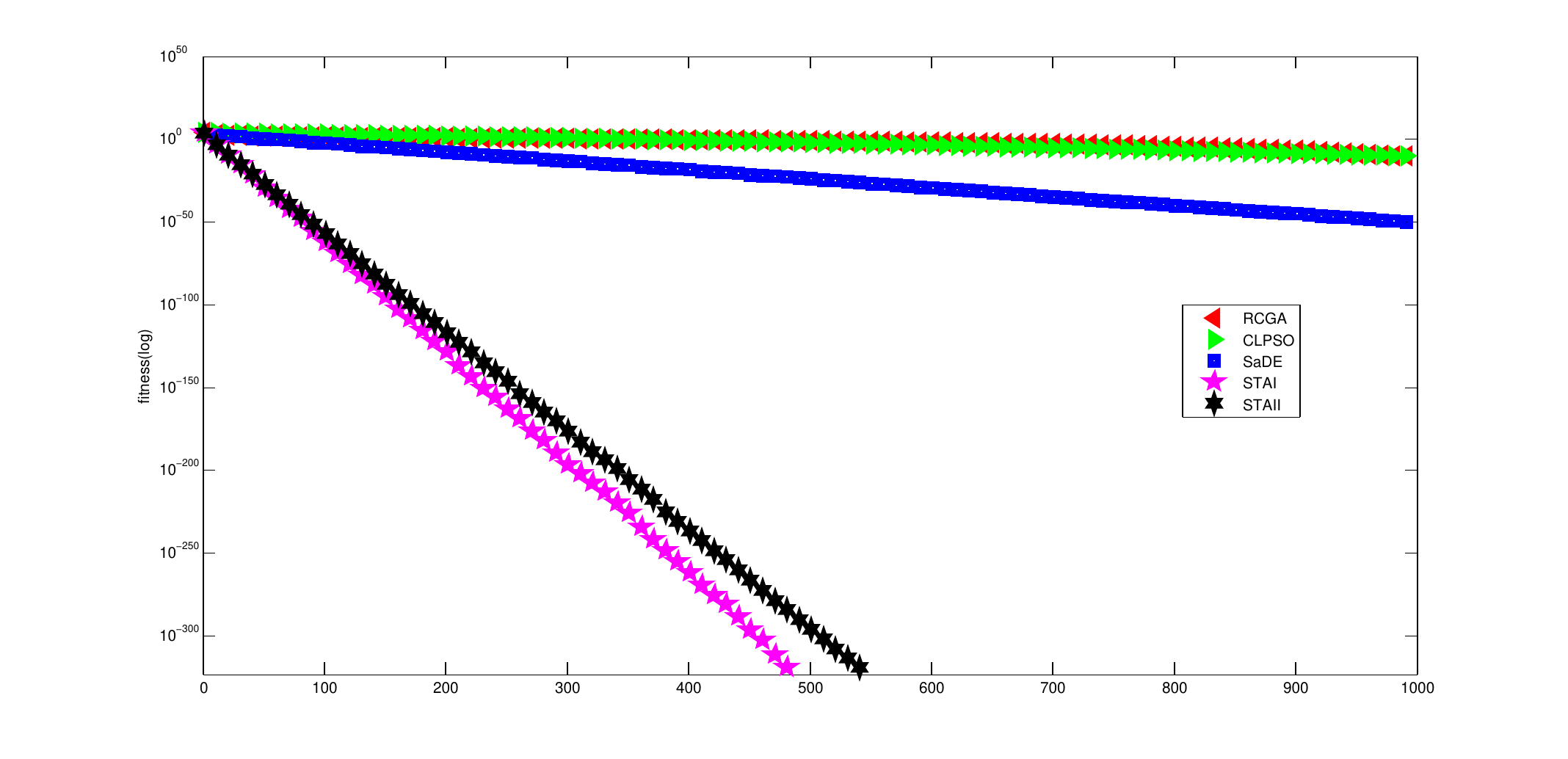}}
\subfloat[$f_2$]{\includegraphics[width=4cm,height=4cm]{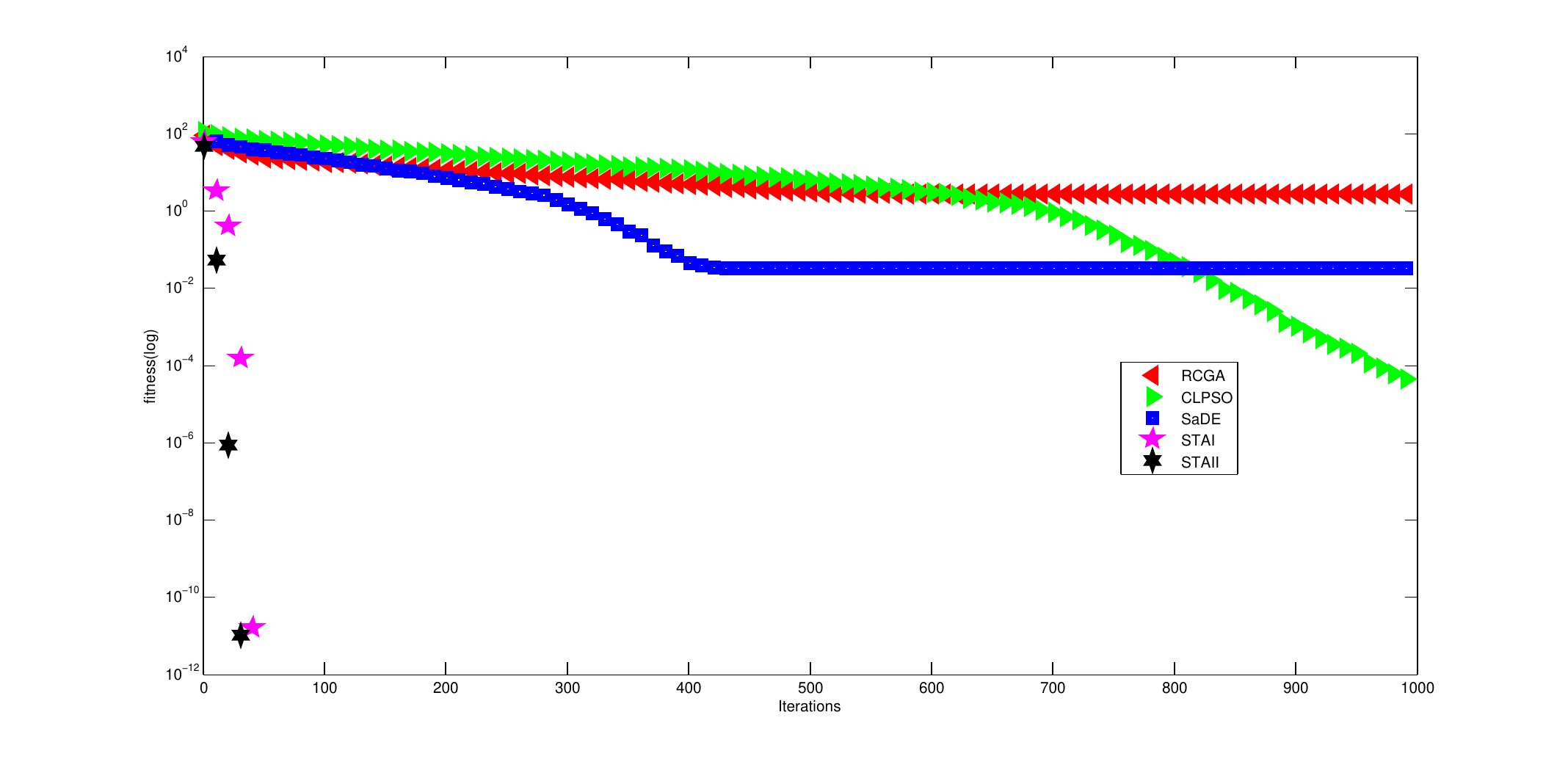}}
\subfloat[$f_3$]{\includegraphics[width=4cm,height=4cm]{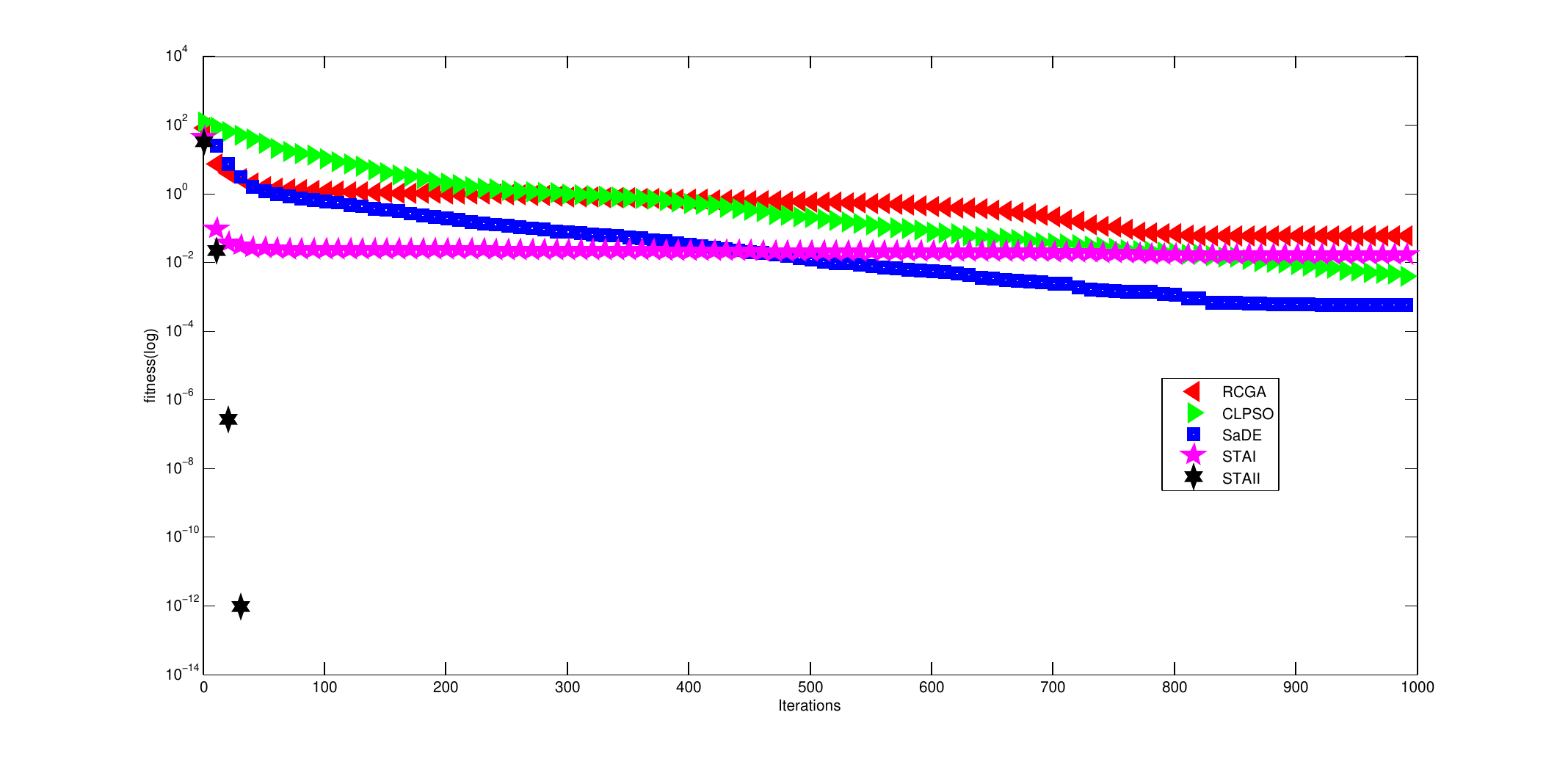}}\\
\subfloat[$f_4$]{\includegraphics[width=4cm,height=4cm]{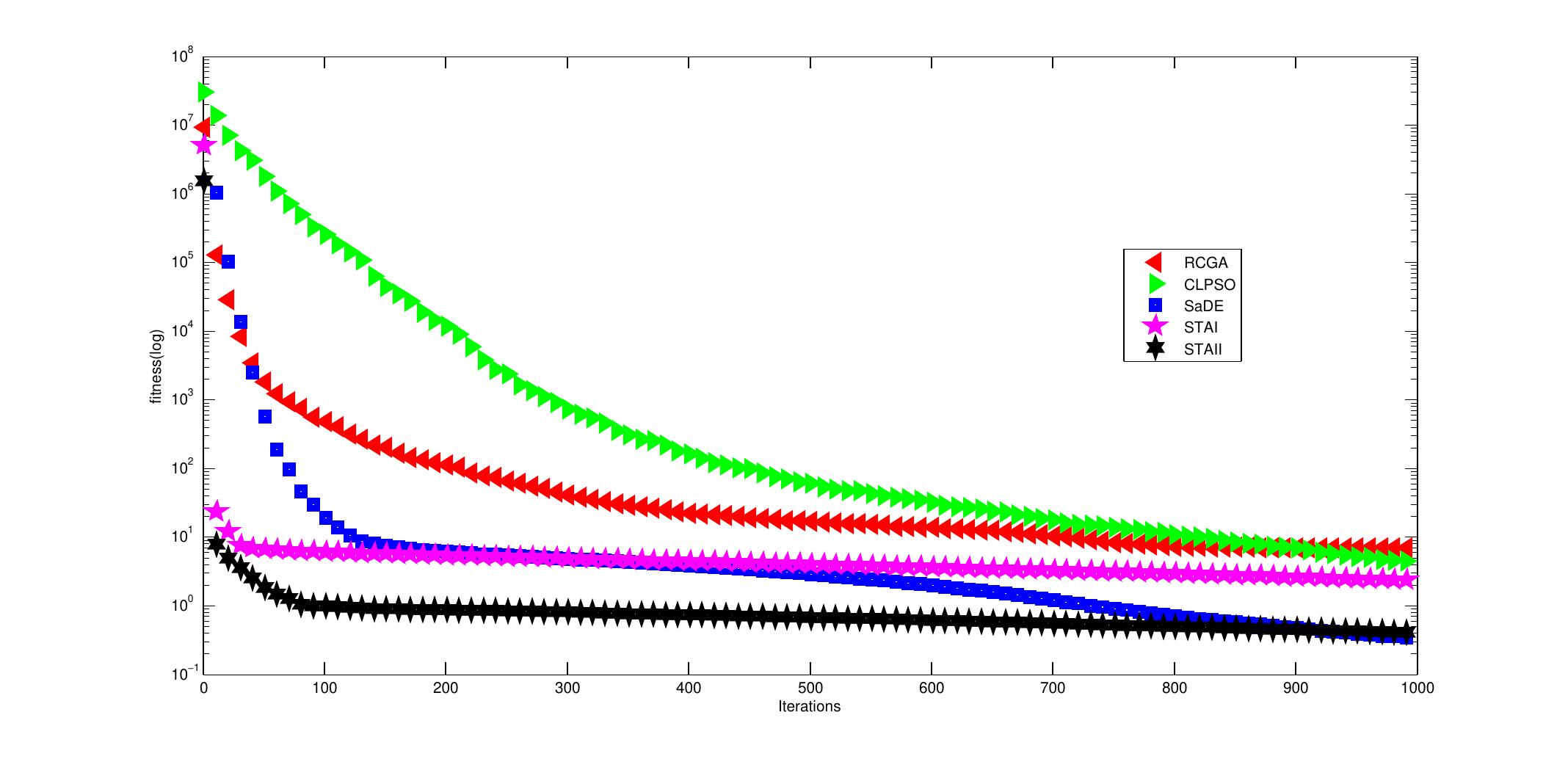}}
\subfloat[$f_5$]{\includegraphics[width=4cm,height=4cm]{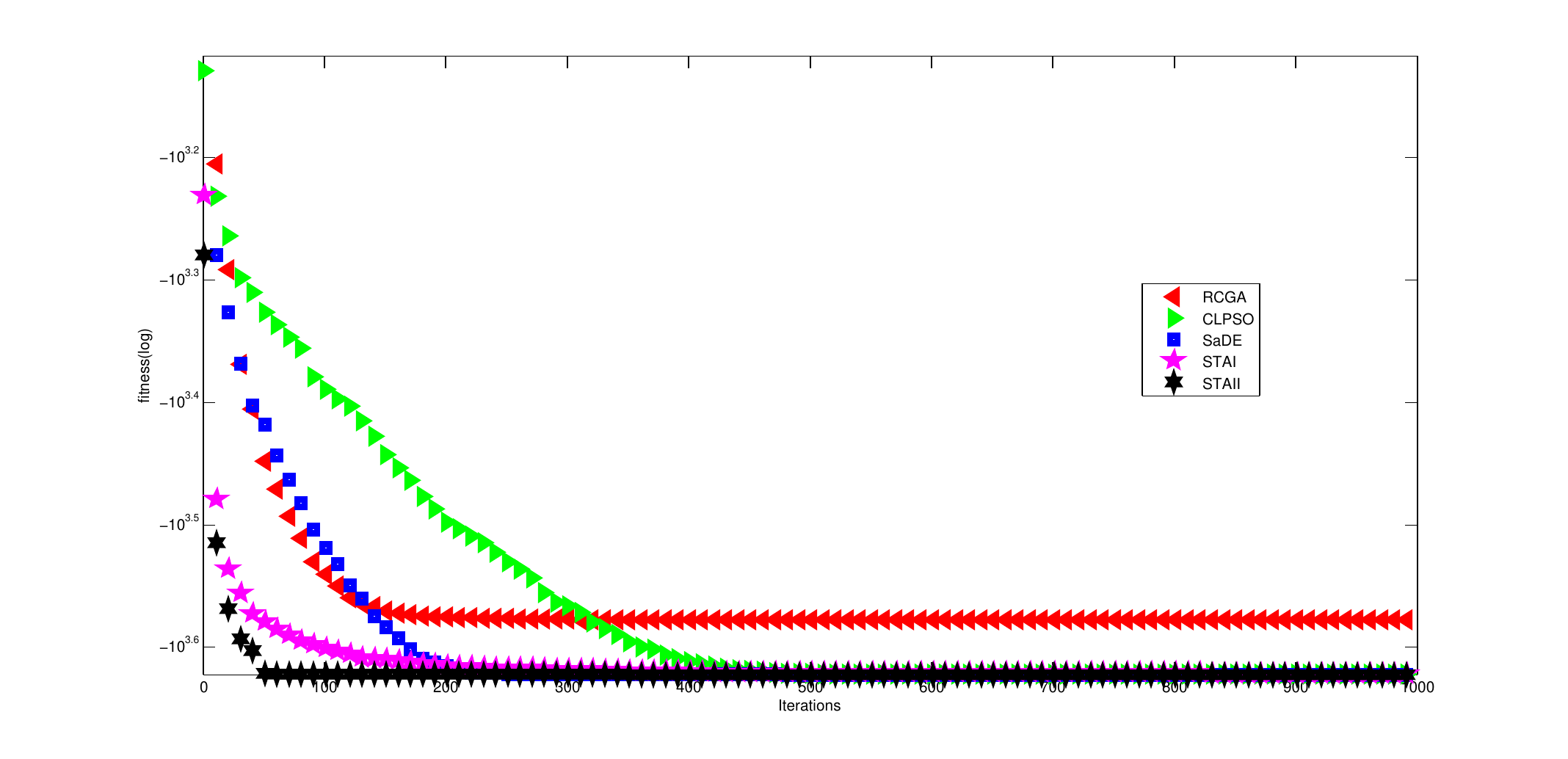}}
\subfloat[$f_6$]{\includegraphics[width=4cm,height=4cm]{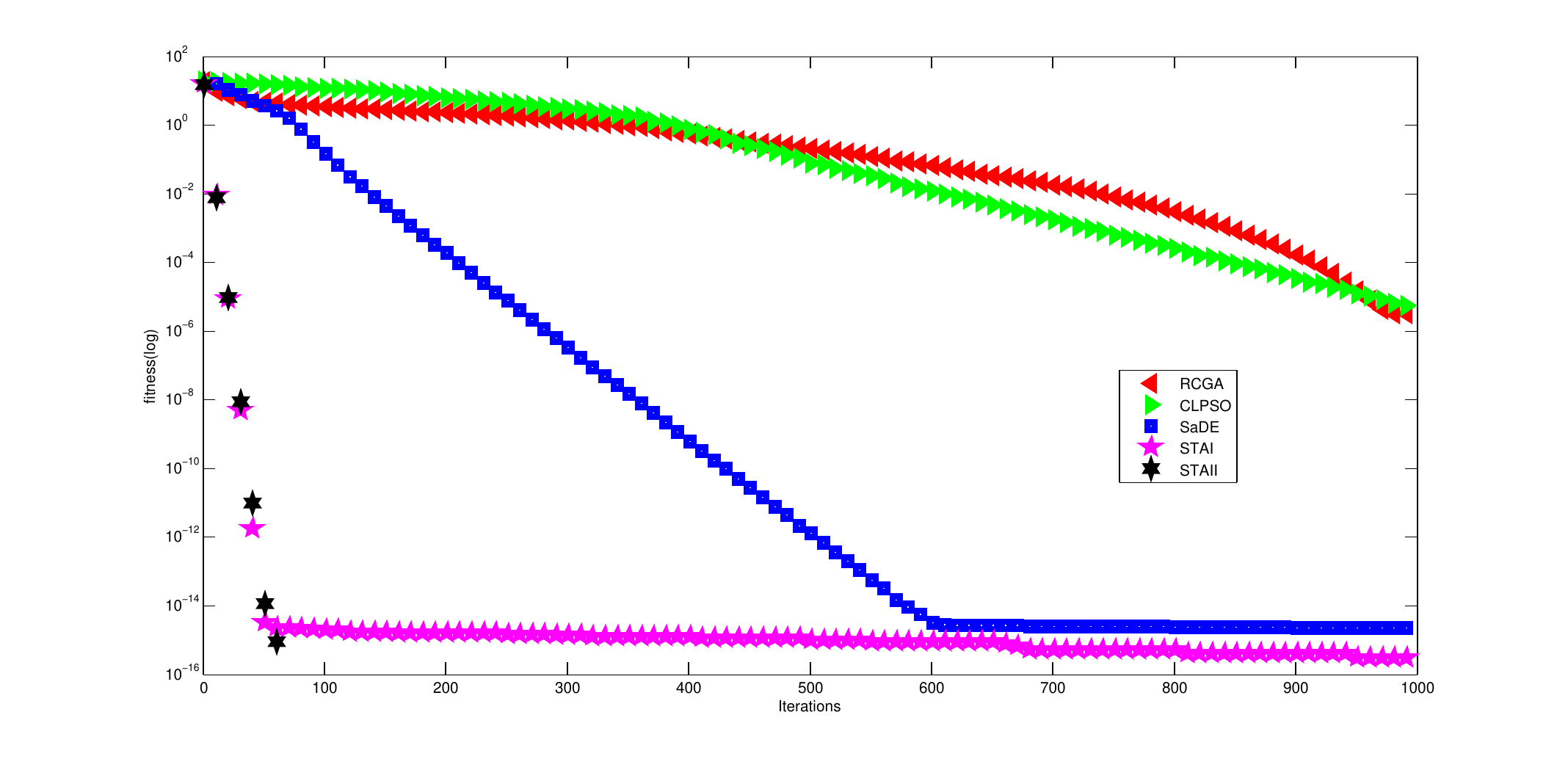}}\\
\subfloat[$f_7$]{\includegraphics[width=4cm,height=4cm]{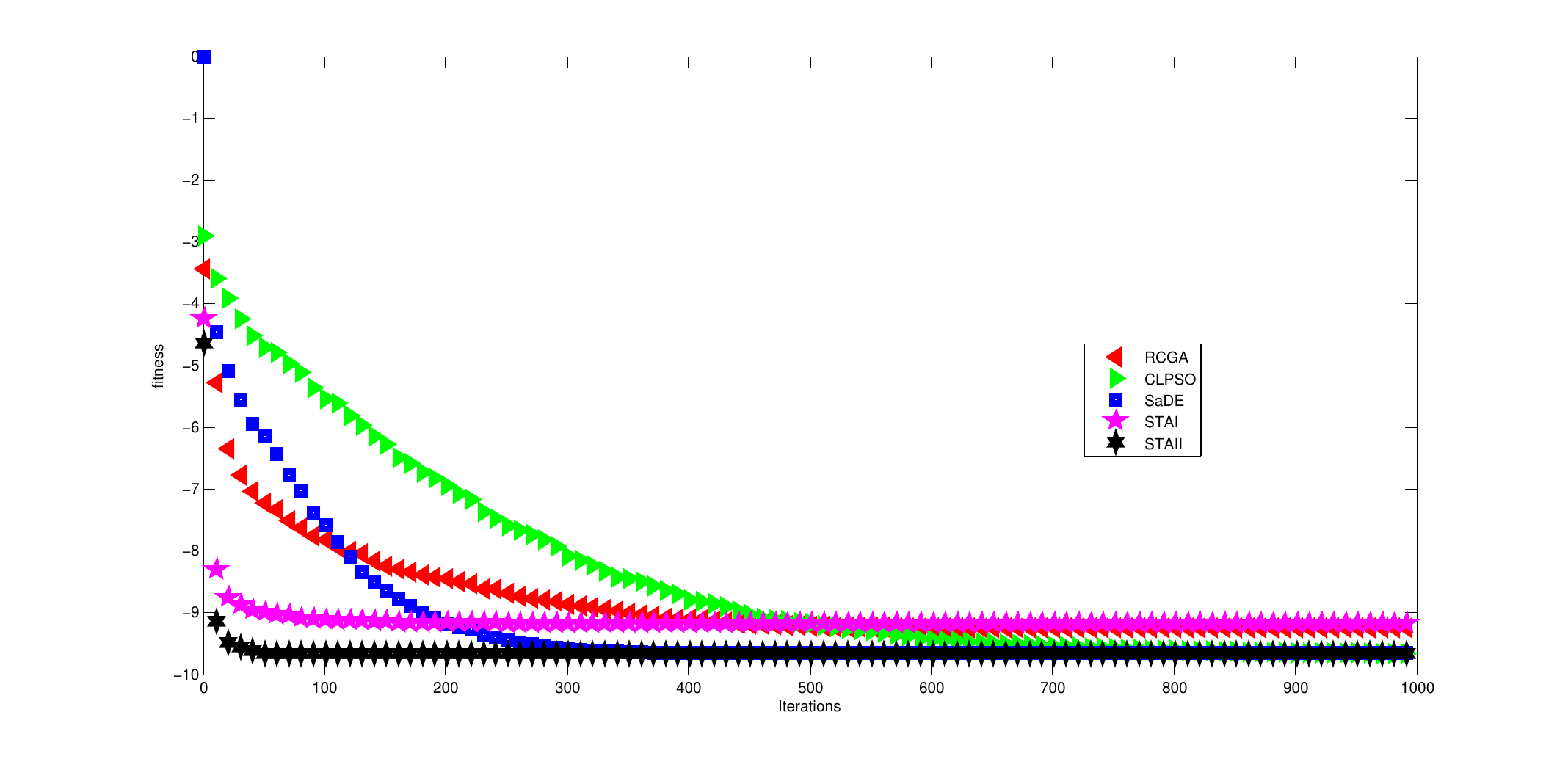}}
\caption{Average fitness of the ten-dimensional functions from $f_1$ to $f_{7}$}
\end{figure}
\\
\indent \textbf{Spherical Function:} as can be seen from the results, all of the algorithms can find the global optimum with high solution precision and have good reliability as well as stability for this function in terms of two and ten dimensions. But the STAI and STAII are able to search much deeper than other three algorithms, which can also be observed in subfigure (A) of Figure 3 and Figure 4. In the subfigure (A), we can see that STAs can converge much faster than the remained methods. While for STAI and STAII, it is found that STAI has a little faster convergence performance than that of STAII.\\
\indent \textbf{Rastrigin Function:} we can see from the results that all of the algorithms can find the global optimum and have
good reliability as well as stability in two dimension. For the ten dimensional problem, the global optimum can also be found by all algorithms; however, STAI
and STAII have better statistical performances than other three algorithms especially described by the \textit{worst}. RCGA and SaDE can not achieve the best occasionally, and the \textit{mean} of RCGA is not satisfactory. From subfigure (B) of Figure 3 and Figure 4, we can also find that STAs converge much faster than other algorithms, and higher solution precision can be obtained. In this time, the process of STAII is slightly faster than that of STAI.\\
\indent \textbf{Griewank  Function:} from the results, we can find that most algorithms have the ability to achieve the best and have both reliability and stability in the two dimensional function except the RCGA. While for the ten-dimension function, these methods are able to find the global optimum but the statistical performances are not satisfactory except STAII, the results of which are excellent. In subfigure (C) of Figure 3 and Figure 4, it can be found that STAII converge fastest and have highest solution precision of all. \\
\indent \textbf{Rosenbrock Function:} all of the algorithms have no problem to find the global optimum for this function in two-dimensional space,
but the \textit{worst} of RCGA indicate that it is not reliable and a bit deficient for the function. Regarding corresponding ten dimensional problem,
only SaDE and STAs can find the best with a low probability. In this case, SaDE achieves best results, followed by STAII. From subfigure (D) of Figure 3 and Figure 4, we can find that STAs still converge faster than other algorithm but with not higher solution precision than SaDE. Compared with STAI, STAII have much better statistical performances, which are indicated by the \textit{worst} and the \textit{mean}.\\
\indent \textbf{Schewefel Function:} as for the function, only RCGA can not find the global optimum for both two and ten dimensions; furthermore, the \textit{median} and \textit{st.dev.} also show that the RCGA is not stable and reliable for this function. Other algorithms achieve the best as well as good reliability and stability because the \textit{st.dev} approaches zero for these methods. From subfigure (E) of Figure 3 and Figure 4, the faster convergence speed belongs to the STAs as well. While for STAs, it shows that STAII converge faster than STAI for the function.\\
\indent \textbf{Ackley  Function:} it seems that all of the algorithms have no problem in finding the global optimum for the function in terms of two and ten dimensions. The statistical performances of results are satisfactory for all methods because the \textit{st.dev.} approaches zero. In subfigure (F) of Figure 3 and Figure 4, we can find that STAs also have faster convergence speed than other algorithms and the solution precision is also higher for STAs when compared with others.\\
\indent \textbf{Michalewicz Function:} all of the algorithms can achieve the global optimum for this function in two dimension, and the statistical performances
is satisfactory in this case. While for the function with ten dimension, only STAII are able to achieve the same statistics as the results in two dimensional function. More specifically, RCGA and CLPSO are not able to find the best. The \textit{worst} of STAI show that it is not reliable sometimes.\\
\indent \textbf{Schaffer Function:} the global optimum can be found by all the algorithms; however, only STAs can achieve reliable and stable performance for this function, as indicated by the \textit{mean} and the \textit{st.dev.}. Other methods fail to find the best occasionally, which is described by the \textit{worst}, that is to say, other algorithm are not reliable for the function. In the case, STAII converge much faster than STAI, as illustrated by subfigure (H).\\
\indent \textbf{Easom Function:} all of the algorithms are able to find the global optimum with a high probability. The \textit{st.dev.} indicates
that the statistical performance are also fine for all methods. The subfigure (I) of Figure 3 shows that STAs have better convergence performance again.\\
\indent \textbf{Goldstein-Price Function:} as described in Table 4, global optimum can be found by all algorithms, the results of which are satisfactory because
the \textit{st.dev.} approaches zero. The subfigure (J) of Figure 3 shows that the convergence speed is fine for all methods but the STAs are much better to some extent.\\
\indent Over all, some explanations can be given on the behavior of average fitness curves. As shown in Figure 3 and Figure 4, the curves of STAs change steadily during the iteration process in most cases, There are two reasons that account for the phenomenon. Firstly, the rotation guarantees the steady decrease of the curves because the rotation factor changes from a maximum value to a minimum value in a periodical way, which prevents current best state from changing sharply. If other transformations do not work, then rotation will help searching in depth with a high precision. Secondly, expansion and translation are beneficial for searching in a new area, while the axesion is proposed to strength the single dimensional search, which are all advantageous for the decrease of the curves.\\
\indent But every once in a while, especially described by the average fitness of $f_4$ and $f_7$ in ten dimension, STAs fail to find the global optimum. As declared in Part 2.2, rotation transformation is used for local search, expansion, translation, and axesion are helpful for global search. In current STAs, their control parameters(transformation factors) are determined by experimental experience for simplicity. The failure of STAs for $f_4$ and $f_7$ occasionally indicate that the global search transformations need to be deeply studied. Anyway, the smaller rotation factor will facilitate the exploitation and the bigger expansion, translation and axesion factors will benefit the exploration, though how to balance them are still pending. Regarding the influence of the \textit{CF}, we can find that STAII has stronger search ability than STAI as the introducing of intermittent exchange. As illustrated by the average fitness curves, the fitness by STAII can still decrease even if that of STAI is already steady, that is to say, the communication strategy can help share information and prevent premature convergence. If the \textit{CF} is larger, more information will be shared, and if the \textit{CF} is small, self development will be enhanced.\\
\indent By the way, the searching time required for STAs is infinity in theory, which is the consequence of random search methodology. However, in practice, we can stop the iteration process by presetting some criteria, for example, the prescribed maximum iterations, or when the fitness is unchanged for a number of times. In this paper, the maximum iterations is used.
\section{Conclusion}
Based on state and state transition, the STA, not only has a simple form but also possess clear geometric significance, which is easy for understanding. Concerning the continuous function optimization problems, it presents the state transformations including rotation transformation, translation transformation and expansion transformation as well as axesion transformation. The paper focuses on the unconstrained optimization problems, and it studies mainly on the approaches of transformations. Furthermore, to enhance its performance in high dimensional functions optimization, communication strategy has been introduced, and the intermittent exchange is proposed to strength the search ability as well as prevent premature convergence. Using 10 benchmark functions for testing, compared with some distinguished optimization algorithms, it shows that STAs have fine performance in terms of global search ability and convergence accuracy, which confirms the effectiveness of the proposed algorithms.\\
\indent On the other hand, distinguished from other population-based algorithms, STA is not originated from simulating natural intelligence, but it takes advantages of the space structure of a function, which opens a new window for optimization. In the paper, control parameters of STAs are not studied deeply, and they are only determined by the experimental experience or for simplicity. In our future work, these problems will be focused on to better develop the state transition algorithms.
\section*{Acknowledgments}
We would like to thank the anonymous referees for their valuable comments and suggestions towards the improvement of this paper.


\end{document}